\documentclass[letterpaper,11pt]{article} 
\usepackage{amsfonts,color,latexsym,amsthm,amssymb,amsmath,epsfig,rotating}

\newtheorem{dfn}{Definition}[section]
\newtheorem{rem}[dfn]{Remark} 
\newtheorem{prop}[dfn]{Proposition}
\newtheorem{thm}{Theorem}[section]
\newtheorem{thma}[dfn]{Theorem} 
\newtheorem{lemma}[dfn]{Lemma}
\newtheorem{cor}[thm]{Corollary}
\newtheorem{ex}[dfn]{Example}
\definecolor{dblue}{rgb}{0,0,.3}

\newlength{\smale}
\newlength{\jmr}
\newlength{\khov}
\newlength{\bernd}
\newcommand{\thth}{{\text{\underline{th}}}}

\newcommand{\np}{{\mathbf{NP}}}

\newcommand{\eps}{\varepsilon}
\newcommand{\Pro}{{\mathbb{P}}}

\newcommand{\Q}{\mathbb{Q}}
\newcommand{\R}{\mathbb{R}}
\newcommand{\C}{\mathbb{C}}
\newcommand{\N}{\mathbb{N}}
\newcommand{\Z}{\mathbb{Z}}
\newcommand{\bO}{\mathbf{O}}

\newcommand{\Zn}{\Z^n}

\newcommand{\Rn}{\R^n}

\newcommand{\vp}{{\psi}}

\newcommand{\Cn}{\C^n}

\newcommand{\Cs}{\C^*}
\newcommand{\Rs}{\R^*}
\newcommand{\cA}{{\mathcal{A}}}
\newcommand{\cH}{{\mathcal{H}}}
\newcommand{\cL}{{\mathcal{L}}}

\newcommand{\cT}{{\mathcal{T}}}

\newcommand{\Rsn}{{(\R^*)}^n}
\newcommand{\Csn}{{(\C^*)}^n}

\newcommand{\bla}{{\overline{\nabla}_{\cA}}} 
\renewcommand{\qed}{$\blacksquare$}
\newcommand{\dia}{$\diamond$}

\newcommand{\cay}{\mathrm{Cay}}
\newcommand{\newt}{\mathrm{Newt}}
\newcommand{\supp}{\mathrm{Supp}}
\newcommand{\conv}{\mathrm{Conv}}

\newcommand{\codim}{\mathrm{codim }}

\begin{document}
\title{\mbox{}\\
\vspace{-1in}
Extremal Real Algebraic Geometry and $\cA$-Discriminants$^1$ } 

\author{
Alicia Dickenstein\thanks{ Departamento de Matem\'atica, FCEN, 
Universidad de Buenos Aires, 
Cuidad Universitaria--Pabell\'on I, 1428 Buenos Aires, Argentina, 
{\tt alidick@dm.uba.ar} \ . Partially supported by UBACYT X042, 
CONICET PIP 5617, and ANPCYT PICT 20569, Argentina.} \and  
J.\ Maurice Rojas\thanks{
Department of Mathematics,  
Texas A\&M University
TAMU 3368, 
College Station, Texas \ 77843-3368,  
USA,  {\tt rojas@math.tamu.edu} \ ,  
{\tt www.math.tamu.edu/\~{}rojas} \ .  
Partially supported by NSF CAREER grant DMS-0349309, NSF REU grant  
DMS-0243822, NSF individual grant DMS-0211458, and Sandia National  
Laboratories.}
\and 
Korben Rusek\thanks{ Texas A\&M University, Department of Mathematics, 
TAMU 3368, College Station, TX \ 77843-3368, {\tt korben@rusek.org} \ . 
Partially supported by the NSF REU program through grant DMS-0243822. } \and 
Justin Shih\thanks{ UCLA Mathematics Department, Los Angeles, CA \ 90095-1555, 
{\tt justin.shih@gmail.com} 
\ . Partially supported by the NSF REU program through grant DMS-0243822. }  } 

\date{\today} 

\maketitle

\mbox{}\hfill
{\em For Askold Georgevich, with admiration and respect.  Happy 60!}  
\hfill\mbox{}


\begin{abstract} 
We present a new, far simpler family of counter-examples to 
Kushnirenko's Conjecture. Along the way, we illustrate  
a computer-assisted approach to finding sparse polynomial  
systems with maximally many real roots, thus shedding light on 
the nature of optimal upper bounds in real fewnomial theory. 
We use a powerful recent formula for the 
$\cA$-discriminant, and give new bounds on the topology of certain 
$\cA$-discriminant varieties. A consequence of the latter result 
is a new upper bound on the number of topological types of certain real 
algebraic sets defined by sparse polynomial equations. 
\end{abstract}

\section{Introduction} 
The algorithmic study of real solutions of systems of multivariate polynomial 
equations is central in science and engineering, as well as in mathematics. 
For instance, entire fields such as Computer Aided Geometric Design and 
Control Theory essentially revolve on basic but highly non-trivial 
questions involving certain structured polynomial systems (see, e.g., 
\cite{rimas,pole}). Furthermore, 
polynomial systems whose real roots lie outside the reach of current 
algorithmic techniques regularly occur in a myriad of industrial 
problems, and many of these problems involve {\bf sparse} polynomial 
systems, i.e., equations with ``few'' monomial terms. Understanding 
the number of real solutions of sparse polynomial equations is 
thus fundamentally important in real algebraic geometry. 

Here we shed light on the difficulty behind 
determining the maximal number of real roots of polynomial 
systems with a fixed number of exponent vectors. We give new, 
dramatically simpler counter-examples (in Theorem~\ref{thm:new} below) to an 
earlier conjectural upper bound of Kushnirenko. A consequence of our 
investigation 
is a precise quantitative statement that (for many {\bf fixed} choices of 
exponent vectors) sparse polynomial systems with maximally many real 
roots are very rare: they lie in 
extremely small chambers --- determined by a suitable discriminant 
variety --- in the space of coefficients (Theorem~\ref{thm:prob} below).  
Moreover, we prove an explicit upper bound on the number of such chambers
(see the proof of Theorem~\ref{thm:disc} below). This in turn implies a new 
upper bound on the number of smooth topological types attainable in families 
of real algebraic sets defined by certain sparse polynomials 
(Theorem~\ref{thm:disc} and Corollary~\ref{cor:sys} below). 

The techniques of our paper actually extend to general $\cA$-discriminants 
and counting topological types of real zero sets of general 
sparse polynomials. The latter results will appear in a forthcoming paper. 
However, the special cases covered here already 
yield new results on extremal real algebraic geometry, which we now 
review in detail. 

\subsection{Background on Extremal Estimates} In a book published in June 
1637, Ren\'e Descartes 
stated that any real univariate polynomial with exactly $m$ monomial terms has 
at most $m-1$ positive roots \cite{descartes}. Unlike the 
behavior of complex roots, Descartes' bound on the number of real roots is 
completely independent of the degree of the polynomial. However, nearly four 
centuries later, we still lack a definitive analogue for systems of 
multivariate polynomials. Great progress was made by Khovanskii and 
Sevastyanov \cite{kho,few} around the late 1970's, 
culminating in an explicit upper bound for the number of 
non-degenerate positive roots of general sparse polynomial systems. This 
bound ---  a very special case of {\bf Khovanskii's Theorem on Real 
Fewnomials} \cite{few} --- revealed that the maximal number of isolated real 
roots of polynomial systems with a fixed number of exponent vectors is 
independent  of the sizes of the exponent vectors. Khovanskii's theory has 
since enabled important advances in many different areas, e.g., 
Hilbert's 16$^\thth$ Problem \cite{kaloshin}, algorithmic complexity  
\cite{gv,butterfly,thresh}, the study of torsion points on algebraic 
curves \cite{cohenzannier}, and model theory \cite{wilkie}, to name but a few 
(see also the conclusion and bibliography of \cite{few}). 

Khovanskii's original bound is now known to be far from tight (see, e.g., 
\cite{tri,bs}). More to the point, finding general {\bf optimal} bounds is 
a decades-old problem whose solution would have significant impact outside, 
as well as inside, real algebraic geometry. Unfortunately, finding optimal 
bounds even for two equations in two unknowns --- 
with just three terms each --- turned out to be difficult 
enough to take close to 20 years to do. 

To clarify this difficulty, consider the following polynomial system,  
which we will call the {\bf Haas system with parameters} $\pmb{(a,b,d)}$:\\ 
\mbox{}\hfill $H_{(a,b,d)}:=\left\{ \begin{matrix} h_1(x,y):= x^{2d}+ay^d-y\\ 
 h_2(x,y):= y^{2d}+bx^d-x\end{matrix} \right.$\hfill\mbox{}\\ 
Letting $\Rs\!:=\!\R\setminus\{0\}$ and letting $\R_+$ denote the 
positive ray, we call the collection  of systems
$\{H_{(a,b,d)}\}_{(a,b,d)\in (\Rs)^2\times\N}$ the {\bf Haas family}.  
\begin{dfn}
Given any $f\!\in\!\R[x_1,\ldots,x_n]$ with
$f(x)\!=\!\sum^m_{i=1}c_ix^{a_i}$ (where the $a_i$ are pair-wise distinct 
and all $c_i$ are nonzero)\footnote{For any $\alpha\!\in\!\Rn$, 
the notation $\alpha\!:=\!(\alpha_1,\ldots,\alpha_n)$ and 
$x^{\alpha}\!:=\!x^{\alpha_1}_1\cdots x^{\alpha_n}_n$ will be understood.}  
we call $f$ a {\bf (real) $n$-variate $m$-nomial}. Also, given 
$f_1,\ldots,f_k$ 
with $f_i$ an $n$-variate $m_i$-nomial for all $i$, we call 
$F\!:=\!(f_1,\ldots,f_k)$ 
a $k\times n$ {\bf fewnomial system of type } $(m_1,\ldots,m_k)$. 
Finally, if the total number of {\bf distinct} exponent vectors among 
the $f_i$ is $m$, then we can also call $F$ an {\bf $m$-sparse $k\times 
n$ fewnomial system}. \dia  
\end{dfn} 

\noindent 
Thus, for example, any system from the Haas family can be referred to as 
(a) a $2\times 2$ fewnomial system of type $(3,3)$, (b) a $2\times 2$ 
trinomial system, or (c) a $6$-sparse $2\times 2$ fewnomial system. 
Note also that $H_{(a,b,d)}$ has the same roots in $(\Rs)^2$ as 
$(h_1(x,y)/y,h_2(x,y)/x)$, which is $5$-sparse. 

The aforementioned special case of Khovanskii's Theorem on Real Fewnomials 
(invoking an improvement observed by Daniel Perrucci \cite{perrucci}) 
states that an $m$-sparse $n\times n$ fewnomial system never has 
more than $(n+1)^{m-1}2^{(m-1)(m-2)/2}$ non-degenerate roots in the 
positive orthant $\Rn_+$. This in turn implies that the maximal number of 
non-degenerate roots in $\R^2_+$ of any $H_{(a,b,d)}$ in the Haas family 
is no more than $5184$, since we can 
replace any $H_{(a,b,d)}$ by a $5$-sparse system with the 
same roots in $\R^2_+$ (cf.\ the preceding paragraph).

Anatoly Kushnirenko, also around the late 1970s, conjectured a significant 
sharpening of Khovanskii's bound: {\bf Kushnirenko's Conjecture} 
was the statement that $n\times n$ fewnomial systems  
of type $(m_1,\ldots,m_n)$ never have more than 
$\prod^{n}_{i=1} (m_i-1)$ non-degenerate roots in the positive orthant 
$\Rn_+$. This conjecture, if true, would have implied that the 
maximal number of non-degenerate roots in $\R^2_+$ of any $H_{(a,b,d)}$ in 
the Haas family is $4$, thus pointing to a 
rather large gap. (It is a simple exercise to 
construct $2\times 2$ trinomial systems having $4$ non-degenerate 
roots in $\R^2_+$.) 

A bit of Gaussian elimination easily reveals that Kushnirenko's conjectural 
bound, if true, would have implied an elegant upper bound of $(m-n)^n$ for the 
number of non-degenerate roots of any $m$-sparse $n\times n$ fewnomial system 
\cite[Prop.\ 1]{tri}. Since $m$-sparse $n\times n$ fewnomial systems 
have {\bf no} isolated roots in $\Rn_+$ when $m\!\leq\!n$ \cite[Prop.\ 1 and 
Thm.\ 4]{tri}, the case where $n$ is fixed and $m$ is large becomes a natural 
question. Kushnirenko's Conjecture (or even an upper
bound of the form $O(m)^{n^{2-\eps}}$ for some $\eps\!>\!0$) 
--- if true --- would have thus been a significant asymptotic improvement to 
Khovanskii's bound. 
\begin{rem}
Curiously, over a different metrically complete field --- the 
$p$-adic rationals $\Q_p$, for any fixed prime $p$ --- it is now known 
that the number of geometrically isolated roots in \linebreak 
$\Q^n_p$ is $O((m-n)\log(m-n))^{3n}$, 
for fixed $n$ and large $m$ \cite{amd}. Also, in the complementary 
setting of fixed $m-n$ and large $n$, the number of non-degenerate roots 
in $\R^n_+$ is now known to be $O(n)^{m-n-1}$ \cite{bs}. \dia 
\end{rem}

According to Bertrand Haas (and conversations of the second author with
Dima Yu.\ Grigor'ev and Askold Khovanskii, on or before September 2000), 
Kushnirenko saw a simple 
counter-example to his conjecture shortly after he stated it in the late 
1970's. Unfortunately, no one ever recorded this counter-example, or the 
identity of its author. Fortunately, Haas proved in 2000 
\cite{haas}, via an ingenious elementary argument, that the system\\ 
\mbox{}\hfill $x^{106}+1.1y^{53}-1.1y$\hfill\mbox{}\\
\mbox{}\hfill $y^{106}+1.1x^{53}-1.1x$,\hfill\mbox{}\\
along with many others with nearby real exponents, has at least $5$ roots in 
$\R^2_+$. Shortly after, Li, Rojas, and Wang proved that {\bf all} 
$2\times 2$ trinomial systems (and, in particular, all systems in the Haas 
family) have at most $5$ isolated roots in $\R^2_+$ \cite{tri}. The latter 
trio of authors also significantly sharpened Khovanskii's bound for certain 
other  families of $n\times n$ sparse polynomial systems. 

Haas' example above was thus the simplest known counter-example to 
Kushnirenko's Conjecture, until the present paper. 
\begin{rem} 
It is interesting to observe that the existence of a pair of 
real bivariate polynomials $F\!:=\!(f_1,f_2)$ --- with $f_1$ a 
trinomial and $f_2$ a tetranomial --- having more than $6$ isolated roots 
in $\R^2_+$, is still an open problem. The maximal number is known to be at 
least $6$, and no larger than $14$ \cite[Thm.\ 1, Assertion (a)]{tri}. \dia
\end{rem}   

\subsection{Main Results} 
\subsubsection{New Counter-Examples and the Probability of Finding One} 
We give a new family of counter-examples far simpler than that of Haas, and 
announce what appear to be many more such families. In particular, while 
Haas found a pair of trinomials of 
degree $106$ (and many more of higher degree), we given an explicit 
cell in $\R^2$ which is naturally identified with an infinite family of 
pairs of trinomials of degree $\pmb{6}$. We have also found 
experimentally $49$ other such cells (consisting of $2\times 2$ systems of 
trinomial of even degree $<106$) which are all counter-examples to 
Kushnirenko's Conjecture as well, but we focus here on the simplest. 
\begin{prop} 
\label{prop:5roots}
Let $E_d\!\subseteq\!\R^2$ denote the set of $(a,b)$ such that
$H_{(a,b,d)}$ has at least $5$ non-degenerate isolated roots in the positive 
quadrant. Then $E_d$ is open and symmetric about the line $\{a=b\}$. \qed  
\end{prop} 
\begin{thm} 
\label{thm:new} 
Following the notation of Proposition~\ref{prop:5roots}, $E_3$ is non-empty 
(and in fact star convex\footnote{Recall that a set 
$S\!\subset\!\Rn$ is {\bf star convex} iff there is a $p\!\in\!S$ such that 
for all $x\!\in\!S$, the (closed) line segment containing $p$ and $x$
is contained in $S$ as well.}). 
In particular, $(\text{\scalebox{.8}[.8]{$\frac{44}{31}$}},
\text{\scalebox{.8}[.8]{$\frac{44}{31}$}})\!\in\!E_3$, and thus the 
$2\times 2$ system $( \ x^6+\frac{44}{31}y^3-y \ , \ 
y^6+\frac{44}{31}x^3-x \ ) $  
has exactly $5$ roots, all non-degenerate, in $\R^2_+$.  
\end{thm} 

\noindent
We point out that while our first verification of our simplest counter-example 
was done via Gr\"obner bases (on the computer algebra system {\tt Maple}), 
we present here a novel (and simpler) numerical verification via 
Smale's {\bf Alpha Theory} \cite{smale,bcss}. 

Our next result reveals that for small $d$, it is highly 
{\bf un}likely that a random choice of $(a,b)$ will yield 
$H_{(a,b,d)}$ as a counter-example to Kushnirenko's Conjecture. 
\begin{thm}
\label{thm:prob} 
Following the notation of Proposition~\ref{prop:5roots}, 
$E_1$ and $E_2$ are empty, and 
$0\!<\!\mathrm{Area}(E_3)\!<\!5.701\times 10^{-7}$. 
In particular, letting $a$ and $b$ be independent identically distributed 
standard real Gaussian random variables, 
$\mathrm{Prob}[(a,b)\!\in\!E_3]\!<\!1.936 \times 10^{-9}$. 
Finally, let $\alpha$ be the smallest positive root of the univariate degree 
$7$ polynomial $C(a)$, which we define as\\
\scalebox{.54}[1]{{\small $823564528378596a^7-7917064766635392a^6
-195134969401159896a^5-668651015982750336a^4-19908809569295316a^3
-564987948607350000a^2+2392425241171875000a-8620460479736328125$}},\\ 
and let $\beta\!<\!\gamma$ be the two 
smallest positive roots of the univariate degree $36$ 
polynomial $A(a)$ in the Appendix. Then the boundary of $E_3$ consists
of exactly $4$ convex arcs, meeting exactly at 4 vertices:\\ 
\mbox{}\hfill$\left\{\sqrt[5]{\alpha}\begin{bmatrix}1\\ 1\end{bmatrix},
\sqrt[5]{\frac{16807}{2916}}\begin{bmatrix}1\\1\end{bmatrix},
\begin{bmatrix}\sqrt[35]{\beta}\\ \sqrt[35]{\gamma}
\end{bmatrix},
\begin{bmatrix}\sqrt[35]{\gamma}\\ \sqrt[35]{\beta}\end{bmatrix} 
\right\}$.\hfill\mbox{}\\ 
To $10$ decimal places, the 
preceding coordinates are\\ 
\mbox{}\hfill$\left\{1.4176759490 \begin{bmatrix}1\\1\end{bmatrix},
1.4195167977\begin{bmatrix}1\\1\end{bmatrix},
\begin{bmatrix}1.4179051055\\1.4182147696\end{bmatrix}
,
\begin{bmatrix}1.4182147696\\1.4179051055\end{bmatrix}
\right\}$.\hfill\mbox{}
\end{thm} 

\noindent 
This paucity of extremal examples can be visualized most easily by 
plotting the regions in the space of coefficients that yield 
$H_{(a,b,d)}$ with a given number of roots in $\R^2_+$.\footnote{The 
fact that systems in the Haas family have no roots on the 
coordinate cross, other than the origin, guarantees constancy 
(if we vary the coefficients while avoiding singularities)
of the number of roots in each quadrant.} 
Below is a sequence of $4$ such plots (for $d\!=\!3$), 
drawn on a logarithmic scale and successively magnified up to a 
factor of about $1700$. \\
\epsfig{file=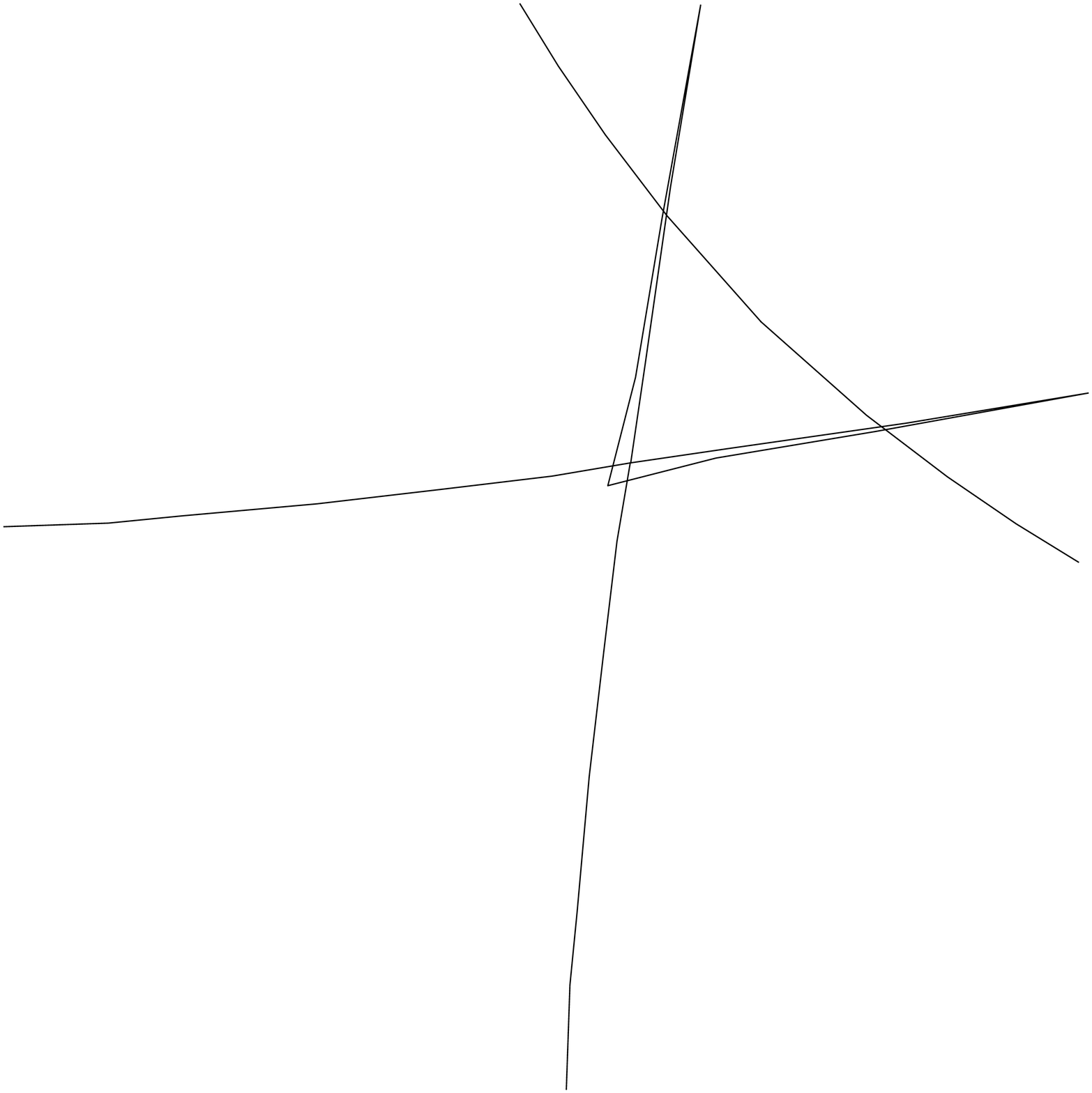,height=1.5in}
\epsfig{file=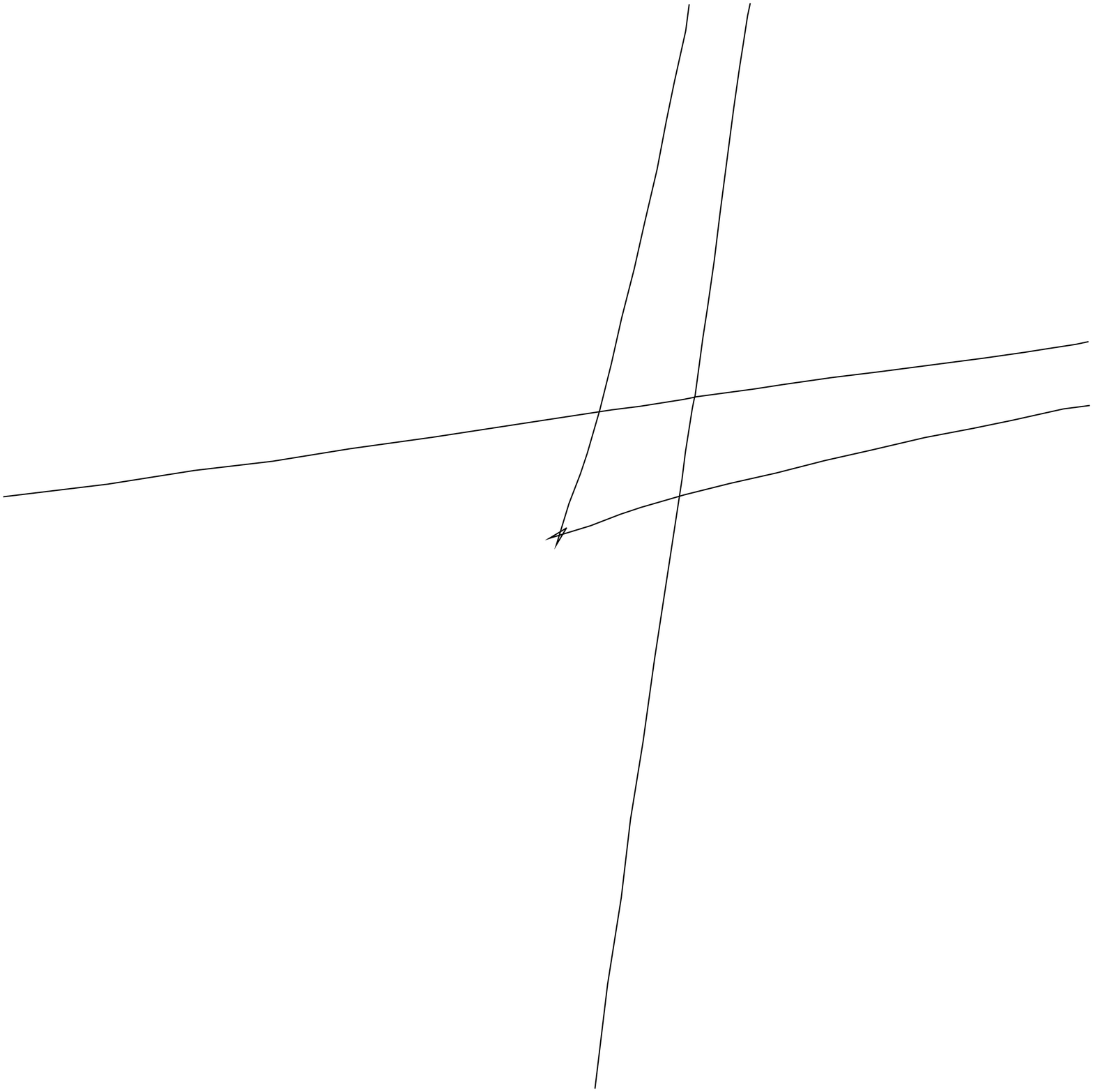,height=1.5in}
\hspace{.75in}\raisebox{.75in}{\epsfig{file=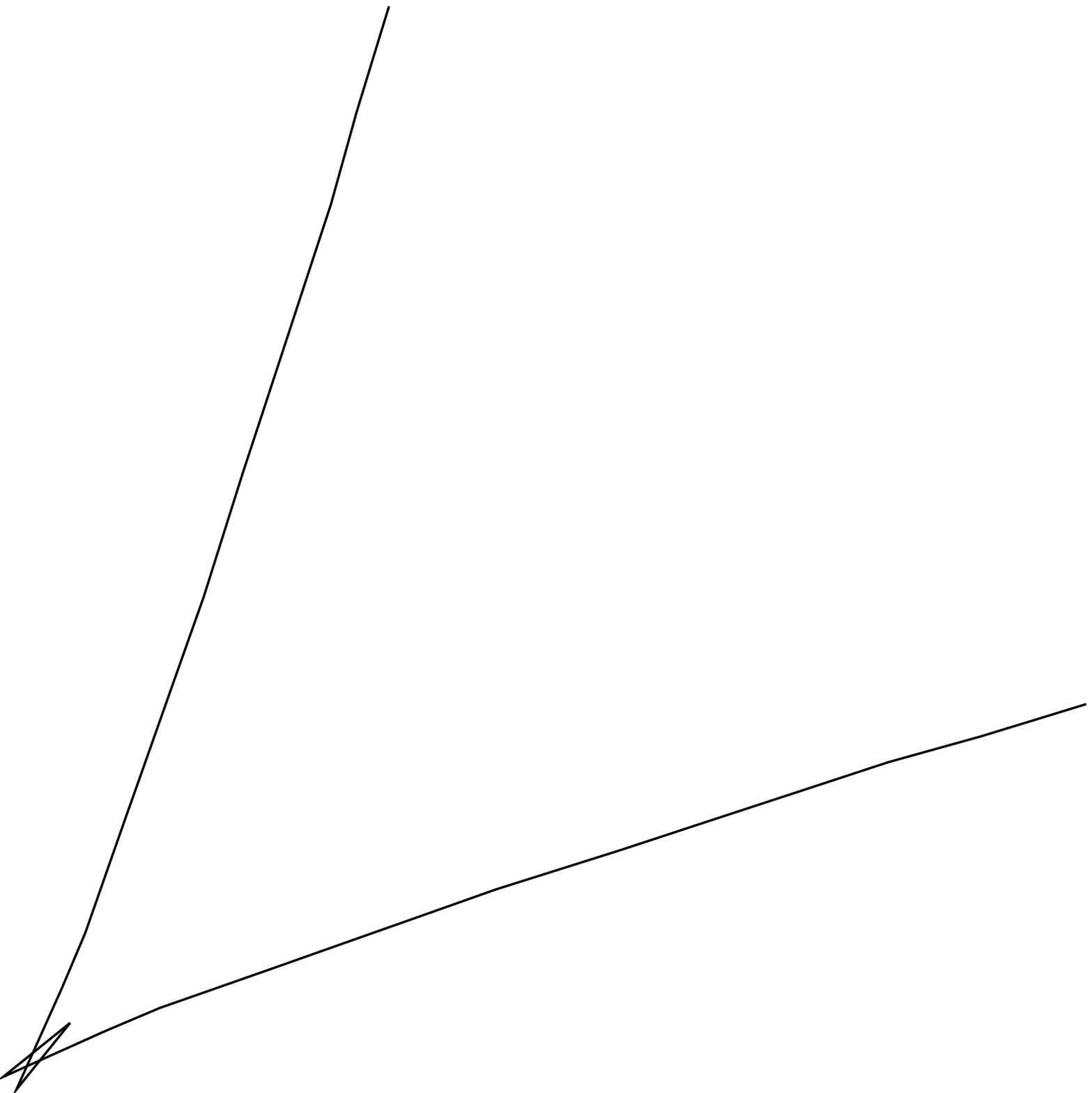,
height=.75in}} 
\epsfig{file=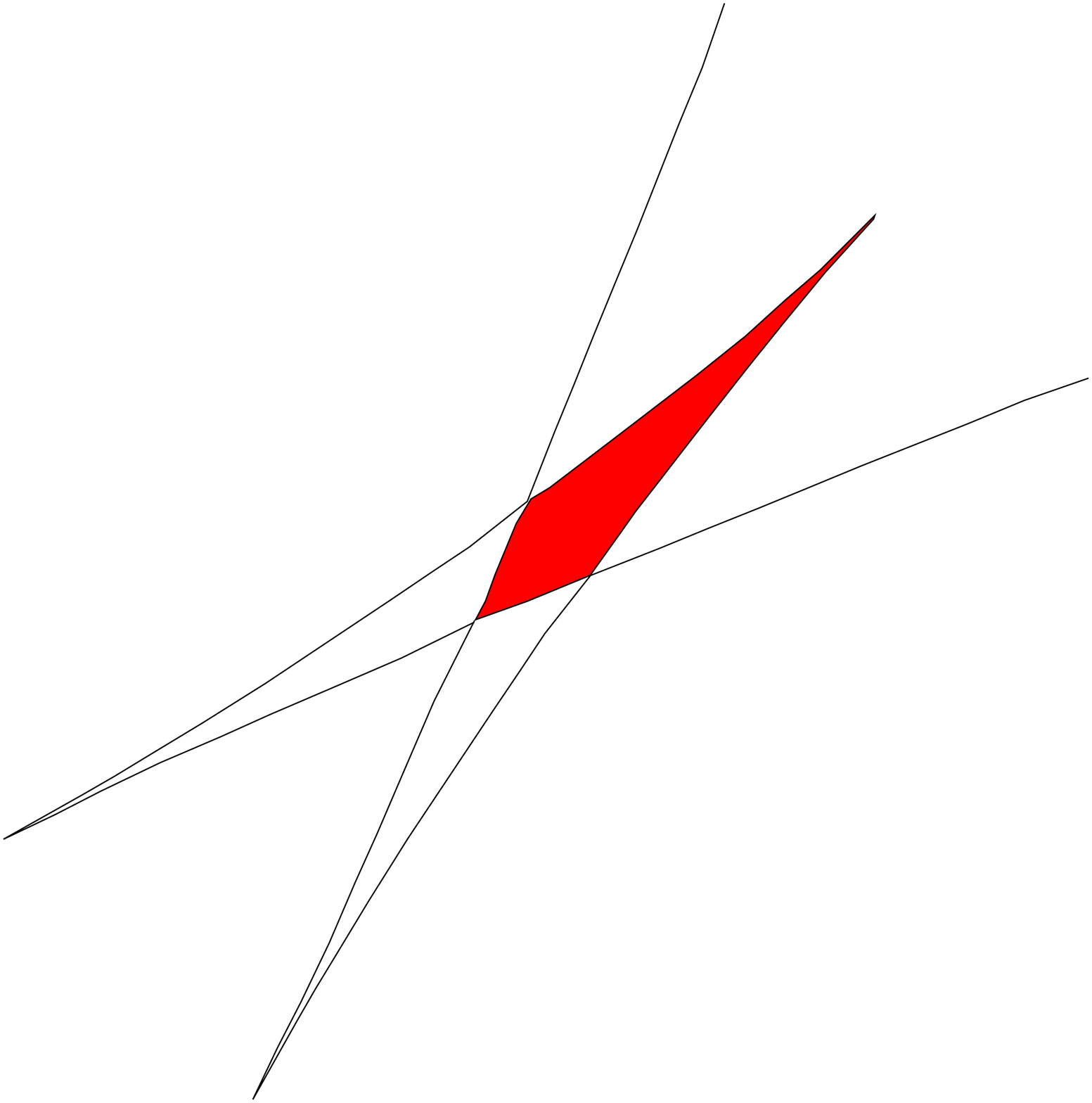,height=1.5in}\\
The small red (or grey) diamond-like region on the right-most 
plot is $E_3$. The connected components of the complement of the underlying 
real curve are regions on which the number of roots of 
$H_{(a,b,3)}$ in $\R^2_+$ --- as a function of $(a,b)$ --- is constant 
(see also the discussion after the proof of Theorem \ref{thm:prob}  
in Section \ref{sec:prob} below). 
In particular, the real curve we see above is the real part of an 
{\bf $\pmb{\cA}$-discriminant variety} --- a tool (reviewed in Section 
\ref{sec:disc} below) underlying our final main result. 

\subsubsection{New Counts for Topological Types and Discriminant Chambers} 
Recall that while a smooth, real, degree $d$ projective plane curve has
at most $1+\text{\scalebox{.7}[.7]{$\begin{pmatrix}d-1\\ 2\end{pmatrix}$}}$
connected components \cite{harnack}, determining the possible nestings of
these ovals --- a piece of the first part of Hilbert's famous 16$^\thth$ Problem
\cite{kaloshin} --- is quite complicated. In more general language, 
this is the determination of possible {\bf diffeotopy types} of 
such curves. 
\begin{dfn} 
Recall that a {\bf diffeotopy} between two sets
$X,Y\!\subseteq\!\Rn$ is a differentiable function $H : [0,1]\times \Rn
\longrightarrow \Rn$
such that $H(t,\cdot)$ is a diffeomorphism for all $t\!\in\![0,1]$,
$H(0,\cdot)$ is the identity on $X$, and $H(1,X)\!=\!Y$. Equivalently,
we simply say that $X$ and $Y$ are {\bf diffeotopic}. \dia 
\end{dfn} 

\noindent 
Note that diffeotopy is a more refined equivalence than diffeomorphism, since
diffeotopy implies an entire continuous family of ``infinitesimal'' 
diffeomorphisms that deform $X$ to $Y$ and back again. Returning to 
nestings of ovals of real degree $d$ projective plane curves, an asymptotic 
formula of $e^{d^2}$ is now known \cite{orevkov}, and the {\bf exact} number is 
currently known (as of late 2006) only for $d\!\leq\!8$. 

Via our techniques here, we can study diffeotopy types in a dramatically 
different setting. 
\begin{dfn}
Given any $n$-variate $m$-nomial $f$, its {\bf support} (or
{\bf spectrum}) --- written $\supp(f)$ --- is its set of exponent vectors.
Also, given any $k\times n$ fewnomial system $F\!=\!(f_1,\ldots,f_k)$,
let $\supp(F)\!:=\!(\supp(f_1),\ldots,\supp(f_k))$.
Finally, we let $Z_+(F)$ (resp.\ $Z^*_\R(F)$) denote the set of roots of $F$ in
$\Rn_+$ (resp.\ $\Rsn$). \dia
\end{dfn}

Given $\cA\!\subset\!\Zn$ with $\#\cA\!=\!n+3$, we let $Q$ denote 
its {\bf convex hull}, and make a mild assumption 
that will be removed in future work (see, e.g., \cite[Sec.\ 3.2]{thresh}): 
To avoid certain technicalities involving topological changes at infinity, we 
assume that the configuration $\cA$ is generic in the sense that the 
intersection of $\cA$ with each facet\footnote{i.e., face of codimension $1$} 
of $Q$ consists of exactly $n$ points. In particular, it is a routine
exercise in polynomial/linear algebra to show that
when the points of $\cA$ are chosen uniformly randomly from
$[-s,s]^n\cap \Zn$, the failure probability of our genericity  
hypothesis decays like $O(1/s^n)$ as $s\longrightarrow+\infty$. 
\begin{thm}
\label{thm:disc}
For any fixed $\cA\!\subset\!\Zn$ 
with $\#\cA\!=\!n+3$ and satisfying the genericity assumption above, 
there are no more than $\left(10+6n+n^2\right)\left(16+8n+n^3
+\frac{e^2+3}{2}\cdot (n+2)^2(n+4)^2\right)$ diffeotopy types for any smooth 
$Z^*_\R(f)$ with $\supp(f)\!=\!\cA$. In particular, the preceding bound 
(a) is no larger than $\frac{26}{5} (n+4)^6$, and  
(b) is completely independent of the  coordinates of $\cA$.  
\end{thm} 

\noindent 
In contrast to the situation for degree $d$ plane curves, it is interesting to 
note that in our $n$-variate $(n+3)$-nomial setting, 
the number of possible diffeotopy types is much 
closer to the maximal number of compact connected 
components of $Z_+(f)$: The latter number was recently shown to be 
no greater than $\left\lfloor\frac{5n+1}{2}\right\rfloor$ \cite{bs}. 
In view of Theorems~\ref{thm:prob} and \ref{thm:disc}, it thus appears that the
intricacies of distinguishing attainable topological types for $Z_+(f)$ ---
for an $f$ with support $\cA\!\subset\!\Zn$ --- might begin
at $\#\cA\!=\!n+3$. 

In particular, for $\#\cA\!=\!n+1$ (resp.\ $\#\cA\!=\!n+2$), there are at most 
$2$ (resp.\ $3$) smooth diffeotopy types for $Z_+(f)$. 
The latter topological bounds are respectively 
proved in Section 3.2 and Theorem 2 of \cite{thresh}, where the resulting 
algorithmic implications are pursued in much greater depth.
We are unaware of any other earlier published explicit bounds --- depending 
solely on $n$ and $\#\cA$ --- for the number of topological types of  
real fewnomial zero sets. That there are any such  
bounds at all is already a non-trivial fact, first observed by 
Lou van den Dries around the 1990s via o-minimality (see, e.g., 
\cite[Prop.\ 3.2, Pg.\ 150]{ominimal}). 
Thierry Zell has informed the authors that results in 
\cite{gabrielov} (on Pfaffian functions and quantifier elimination for 
fewnomials) appear to imply an upper bound of 
$n^{O(n^2)}2^{O(n^4)}$ for the number of 
corresponding smooth topological types in our fewnomial setting 
above.\footnote{ As this paper was being finalized,
Saugata Basu informed the authors that a similar (exponential) bound is likely  
also provable via \cite{basu}, but with the advantage that the 
support can be varied when $n$ and $\#\cA$ are fixed.}  
 
Our polynomial bound above is thus a great improvement. 
It is also likely that our bound above can be 
improved even further: the main example we explore in this 
paper has $n\!=\!3$ and just $15$ representative chambers, even though 
the underlying discriminant has over $58$ monomial terms. 

As a consequence of our last theorem, we get a {\bf constant} 
upper bound for the number of smooth topological types of 
families of real algebraic surfaces defined by ``honest'' hexanomials. 
\begin{cor} 
\label{cor:sys} 
Fix $\cA\!\subset\!\Z^3$, satisfying the assumption preceding 
Theorem \ref{thm:disc} and with cardinality $6$. 
Consider the family of polynomials $g$ with support $\cA$, 
which we will identify with $(\Rs)^6$. 
Then the number of smooth diffeotopy types of 
any such $Z^*_\R(g)$ is no more than $237920$. \qed 
\end{cor} 

\noindent 
\scalebox{.96}[1]{Theorems~\ref{thm:disc}, \ref{thm:new}, and \ref{thm:prob} 
are respectively proved in Sections 
\ref{sec:pfdisc}, \ref{sec:new}, and \ref{sec:prob}. The 
main underlying tools}\linebreak
\scalebox{.94}[1]{are an important recent parametric 
formula for certain $\cA$-discriminant varieties \cite[Prop.\ 4.1]{dfs}, and}
\linebreak  
\scalebox{.92}[1]{recent advances in quantitative estimates for 
{\bf sheared binomial systems} (cf.\ Section \ref{sec:disc}) \cite{tri,bs}.}  

\section{Background on $\cA$-Discriminants} 
\label{sec:disc} 
The standard reference for $\cA$-discriminants is \cite{gkz94}. For our 
purposes, we will modify a few 
notions, present motivating examples, and quote some more recent results as 
well. But first, let us recall the following notation (see \cite{loeser} or 
\cite[Ch.\ 1, 9--11]{gkz94} for further background). 
\begin{dfn} 
Given any $\cA\!=\!\{a_1,\ldots,a_m\}\!\subset\!\Zn$ of cardinality $m$, 
we let $X_\cA$ --- the {\bf (projective) toric 
variety associated to $\cA$} --- be the closure of the parametrized 
subvariety\linebreak  
$\{[t^{a_1}:\cdots:t^{a_m}] \; | \; t\!=\!(t_1,\ldots,t_n)\!\in\!\Csn\}$ of 
$\Pro^{m-1}_\C$.  Finally, we define $\nabla_\cA$ --- the 
{\bf $\cA$-discriminant variety} --- to be the closure of the set of all 
$[c_1:\cdots :c_m]\!\in\!\Pro^{m-1}_\C$ 
such that the hyperplane $\{c_1z_1+\cdots +c_mz_m\!=\!0\}$ intersects 
a regular point of $X_\cA$ with a tangency. \dia 
\end{dfn} 

\noindent 
$\nabla_\cA$ also happens to be the closure of 
those $[c_1:\cdots :c_m]\!\in\!\Pro^{m-1}_\C$ such that the complex zero  
set of $\sum^m_{i=1} c_i x^{a_i}$ possesses a singularity in 
$\Csn$ \cite[Prop.\ 1.1]{gkz94}. 
\begin{ex} 
It is a simple exercise to verify that $\nabla_{\{0,1,2\}}$ is 
the set of all $[a:b:c]\!\in\!\Pro^2_\C$ such that the quadratic 
polynomial 
$f(x)\!:=\!a+bx+cx^2$ has degree $<\!1$ or a double root in $\C$. In 
particular, $\nabla_{\{0,1,2\}}$ is the projectivized complex zero 
set of $b^2-4ac$, and one can check that the affine real zero set of 
$b^2-4ac$ in $\R^3$ is a double cone. Topologically,  
$\Pro^2_\R\cap\nabla_{\{0,1,2\}}$ is thus a circle. \dia 
\end{ex}  

\noindent
In particular, $\nabla_\cA$ is an 
irreducible algebraic variety defined over $\Z$ (see, e.g., 
\cite[Prop.\ 1.3, Pg.\ 15]{gkz94}), and this motivates the following 
important definition. 
\begin{dfn} 
When $\codim \nabla_\cA\!=\!1$, we define (up to sign) the {\bf 
$\cA$-discriminant}, $\Delta_\cA\!\in\!\Z[c_1,\ldots,c_m]$, to be the 
(irreducible) defining polynomial of $\nabla_\cA$. Otherwise, we set 
$\Delta_\cA\!:=\!1$. \dia 
\end{dfn}
\begin{ex}  
Continuing our last example, we can take   
$\Delta_{\{0,1,2\}}\!=\!b^2-4ac$. In particular, observe that for 
{\bf real $(a,b,c)$}, $f$ has exactly $0$, $1$, or $2$ roots 
according as $\Delta_{\{0,1,2\}}(a,b,c)$ is negative, zero, or positive. 
Note also that the classical quadratic formula tells us that the 
roots of $f$ are differentiable functions of the coefficients $(a,b,c)$, 
provided $[a:b:c]\!\not\in\!\nabla_{\{0,1,2\}}$. \dia 
\end{ex} 
\begin{ex} 
For $\cA\!:=\!\{(0,0),(1,0),(2,0),(0,1)\}$, it is easily checked that 
$\nabla_\cA$ is exactly\linebreak 
$\{[a:b:c:d]\; | \; b^2-4ac\!=\!d\!=\!0\}$ and thus 
$\Delta_\cA\!=\!1$. \dia 
\end{ex}

$\cA$-discriminants are central in computational algebraic geometry, 
containing all known multivariate resultants as special cases (see, e.g., 
\cite[Prop.\ 1.7, pg.\ 274]{gkz94}). More to the point, they are 
notoriously difficult to compute: (1) detecting just their 
vanishing is $\np$-hard already for $\cA\!\subset\!\Z^2$ 
\cite{plaisted,koirandim}, 
and (2) $\cA$-discriminants can have many monomial terms, already 
for $\cA\!\subset\!\Z$ and $\#\cA\!=\!4$ \cite[Sec.\ 1.2]{rojasye} (see 
also Example \ref{ex:big} of the next section). 
 
Relevant to our applications, the fact that the vanishing of $\Delta_\cA$ 
determines when certain hypersurfaces possess singularities readily implies 
that the real complement of $\nabla_\cA$ can be used to encode the number of 
real roots of certain 
polynomial systems. In particular, let us call any connected component 
of $\Pro^{\#\cA-1}_\R\!\setminus\!\nabla_\cA$ a(n) 
{\bf ($\cA$-)discriminant chamber}. 
Lemma~\ref{lemma:chamber} of Section \ref{sub:cayley} below (see also 
\cite[Ch.\ 11, Sec.\ 5]{gkz94}) 
relates the number of $\cA$-discriminant chambers to the number
of smooth topological types attainable by the real zero sets of certain 
families of sparse polynomial systems. 
Understanding discriminant chambers is thus a feasible route 
toward understanding the maximal number of real roots of $F$. 
\begin{rem} 
When $\#\cA\!\in\!\{n+1,n+2\}$, $\cA$-discriminant chambers turn out to have
a very simple structure: just one (resp.\ at most two) 
chamber(s) for $\#\cA\!=\!n+1$ (resp.\ $\#\cA\!=\!n+2$) 
\cite[Prop.\ 1.8, Pg.\ 274]{gkz94}. Hence our focus on $\#\cA\!=\!n+3$
throughout this paper. \dia
\end{rem}

\subsection{Studying $\cA$-Discriminant Chambers} 
In what follows, we will always assume that $\cA$ is a subset 
of $\Zn$, of cardinality $m$, that {\bf affinely generates} 
$\Zn$. We will also frequently assume that $\bO\!\in\!\cA$;  
in which case $\cA$ affinely generating $\Zn$ is equivalent 
to $\Zn$ being generated by the set of all integral linear combinations 
of the remaining\linebreak 
\scalebox{1}[1]{vectors in $\cA$. In general, these assumptions  
can easily be enforced by merely translating $\cA$ and, if}\linebreak 
\scalebox{.95}[1]{necessary, applying the {\bf Hermite factorization} 
\cite{storjo} for integral matrices. 
$\cA$-discriminants happen}\linebreak 
to be invariant (modulo a permutation of 
coordinates) under affine transformations 
of $\cA$ that are injective and integral \cite[Prop.\ 1.4, Pg.\ 272]{gkz94}. 
For now, let us observe a simple example. 
\begin{ex}
Note that while $\cA\!:=\!\{3,30,57\}$ does {\bf not} affinely generate $\Z$ or 
contain the origin, this $\cA$ is a {\bf shifted multiple} of a set that does 
both: $\{0,1,2\}$. 
In particular, the polynomial $ax^{57}+bx^{30}+cx^3$ has a 
degenerate root in $\Cs$ iff the polynomial $ax^2+bx+c$ has a 
degenerate root in $\Cs$. Thus, we clearly have $\nabla_{\{3,30,57\}}\!=\!
\nabla_{\{0,1,2\}}$ and $\Delta_{\{3,30,57\}}\!=\!\Delta_{\{0,1,2\}}$. \dia 
\end{ex} 

\vspace{-.9cm}
\begin{rem} 
\scalebox{.87}[1]{We will frequently abuse notation by also identifying $\cA$ 
with the 
$n\times m$ matrix \scalebox{1}[.7]{$\begin{bmatrix}a_{11} & \cdots & a_{1m}\\
\vdots & \ddots & \vdots \\ a_{n1} & \cdots & a_{nm}\end{bmatrix}$},}

\vspace{-.4cm} 
\noindent
where $a_i\!=\!\text{\scalebox{1}[.7]
{$\begin{bmatrix}a_{1i}\\ \vdots \\ a_{ni}\end{bmatrix}$}}$
for all $i$. Also, given vectors $v_1,\ldots,v_m\!\in\!\Rn$, 
we will let $[v_1,\ldots,v_m]$ denote the\linebreak 
\scalebox{.97}[1]{$n\times m$ matrix whose 
$j^\thth$ column is $v_j$. Finally, we will use the notation 
$z^{[v_1,\ldots,v_m]}\!:=\!(z^{v_1},\ldots,z^{v_m})$. \dia}  
\end{rem}   

\vspace{-.2cm}
The combinatorics of exponents will thus be particularly important throughout 
our \linebreak 
development. Continuing this theme, we will define notions useful for 
simplifying discriminant chambers. However, before going into further 
definitions, let us first motivate the need for simplification via a more 
intricate $\cA$-discriminant example.

\vspace{-.2cm} 
\begin{ex} 
\label{ex:big} 
Suppose $\cA\!:=\!\text{\scalebox{1}[.7]{$\left\{
\begin{bmatrix} 6 \\ 0 \\ 0 \end{bmatrix},
\begin{bmatrix} 0 \\ 3 \\ 0 \end{bmatrix},
\begin{bmatrix} 0 \\ 1 \\ 0 \end{bmatrix},
\begin{bmatrix} 0 \\ 6 \\ 1 \end{bmatrix},
\begin{bmatrix} 3 \\ 0 \\ 1\end{bmatrix},
\begin{bmatrix} 1 \\ 0 \\ 1\end{bmatrix}\right\}$}}$.  
The real part of $\nabla_\cA$ is at the heart of Theorem~\ref{thm:new}, and it 
is interesting to observe that $\Delta_\cA$ happens to be rather 
large: $\Delta_\cA(1,a,-1,1,b,-1)$ is (up to sign)...\\
{\scriptsize
$1102507499354148695951786433413508348166942596435546875 
-516440160351044111358464119738658142157348733522052\,{a}^{35}$
\linebreak 
\mbox{}\hfill $ + \text{ {\small 54 additional monomial terms of comparable 
size }} $  
\hfill\mbox{}\\
$ -24519711093887016527058411574716512472434688\,{b}^{46}{a}^{39}+
82754024941868680778822139064668229594467072\,{b}^{47}{a}^{33}$.}\\ 
In particular, $\Delta_{\cA}\left(1,\frac{44}{31},-1,1,\frac{44}{31},-1\right)\!
\neq\!0$, and (via the proof of Lemma \ref{lemma:chamber} below) this is 
equivalent to the fact that $H_{(44/31,44/31,3)}$ has no 
degenerate roots.  In any event, it should be clear that we need a more 
efficient means of addressing $\nabla_\cA$. \dia 
\end{ex} 

A beautiful recent (re)discovery is the fact that while $\Delta_\cA$ can 
be unwieldy, $\nabla_\cA$ always admits a compactly expressible 
parametrization: the {\bf Horn-Kapranov Uniformization} (see also 
\cite{kapranov,passare}). 
\begin{thma} 
\label{thm:dfs} 
\cite[Prop.\ 4.1]{dfs} 
Given $\cA\!:=\!\{a_1,\ldots, a_m\}\!\in\!\Zn$, the discriminant locus 
$\nabla_\cA$ is exactly the closure of\\ 
\mbox{}\hfill
$\left\{ \left[u_1t^{a_1}:\cdots:u_mt^{a_m}\right]\; \left| \; u\!:=\!(u_1,
\ldots,u_m)\!\in\!\C^m,\ 
\cA u \!=\!\bO, \ \sum^m_{i=1}u_i\!=\!0, \ 
t\!=\!(t_1,\ldots,t_n)\!\in\!\Csn\right. \right\}$. \qed\hfill\mbox{}
\end{thma} 
\begin{ex}
\label{ex:big2} 
Continuing Example \ref{ex:big}, it is easily checked that the set of 
vectors $\{(u_1,\ldots,u_6)\}$ needed to form the parametrization from 
Theorem~\ref{thm:dfs} is a vector space with basis\\ 
\mbox{}\hfill$\{(-2, 35, -33 , -12, 0,12),
(-2, 11 , -9 , -4, 4 ,0 )\}$.\hfill\mbox{}\\ 
We thus obtain that, in spite of the 
huge formula for $\Delta_\cA(1,a,-1,1,b,-1)$ we saw earlier, 
$\nabla_\cA\!\subset\!\Pro^5_\C$ is exactly the closure of\\ 
\scalebox{.9}[1]{$\left\{ \left[-(2\lambda+2)t^6_2:(35\lambda+11)t^3_3:
-(33\lambda+9)t_3:-(12\lambda+4)t_1t^6_3:4t_1t^3_2:12\lambda t_1t_2 \right]\; 
\left| \; \lambda\!\in\!\C, \ (t_1,t_2,t_3)\!\in\!(\Cs)^3\right. \right\}$. 
\dia} 
\end{ex} 

Computationally, however, we will need to 
express $\nabla_\cA$ in an even more efficient manner. In particular, 
if we are studying the topology of the zero set of a polynomial as we vary 
its coefficients, we should certainly take advantage of the various 
homogeneities that preserve the topology of the underlying zero set. 
\begin{ex}
\label{ex:big3} 
Returning to Examples \ref{ex:big} and \ref{ex:big2}, note that 
the set of exponent vectors of 
$g(x,y,v)\!:=\!c_1x^6+c_2y^3+c_3y+v(c_4y^6+c_5x^3+c_6x)$ is exactly $\cA$, 
and that the topology of $Z^*_\R(g)$ is preserved 
under nonzero scalings of the coefficient vector $(c_1,\ldots,c_6)\!\in\!
(\Rs)^6$, and nonzero scalings of the variables $v,y,z$. In particular, 
if we would like to find $(\alpha,\beta,\gamma,\delta)\!\in\!(\Rs)^4$ such 
that $\delta g(\alpha x,\beta y,\gamma v)\!=\!x^6+ay^3-y+v(y^6+bx^3-x)$ for 
some real $a$ and $b$, then we must clearly solve the binomial system\\
$(\alpha,\beta,\gamma,\delta)^{\text{\scalebox{.7}[.5]
{$\begin{bmatrix} 
6 & 0 & 0 & 1 \\ 
0 & 1 & 6 & 0 \\
0 & 0 & 1 & 1 \\
1 & 1 & 1 & 1 \end{bmatrix}$}}}\!=\!(c^{-1}_1,-c^{-1}_3,c^{-1}_4,-c^{-1}_6)$. 
By multiplying and dividing equations (mimicking Gaussian elimination), 
one can then derive that\\
\mbox{}\hfill $a\!=\!\frac{c_2}{c_1}\left(\frac{-c_3}{c_1},
\frac{c_4}{c_1},\frac{-c_6}{c_1}\right)^{-\text{\scalebox{.7}[.5]
{$\begin{bmatrix}
-6 & -6 & -5  \\
1 & 6 & 0  \\
0 & 1 & 1  \end{bmatrix}^{-1}
\begin{bmatrix}
-6\\ 3\\ 0
\end{bmatrix}$}}} $ and $b\!=\!  
\frac{c_5}{c_1}\left(\frac{-c_3}{c_1},
\frac{c_4}{c_1},\frac{-c_6}{c_1}\right)^{-\text{\scalebox{.7}[.5]
{$\begin{bmatrix}
-6 & -6 & -5  \\
1 & 6 & 0  \\
0 & 1 & 1  \end{bmatrix}^{-1}
\begin{bmatrix}
-3\\ 0\\ 1  
\end{bmatrix}$}}} 
$.\hfill\mbox{}\\ 
Note in particular that the determinants of the matrices 
\scalebox{.7}[.7]{$\begin{bmatrix}
6 & 0 & 0 & 1 \\
0 & 1 & 6 & 0 \\
0 & 0 & 1 & 1 \\
1 & 1 & 1 & 1 \end{bmatrix}$}  
 and 
\scalebox{.7}[.7]{
$\begin{bmatrix}
-6 & -6 & -5 \\
1 & 6 & 0 \\
0 & 1 & 1 \end{bmatrix}$} 
are both {\bf odd} and differ only in sign. So one can check by hand (or via 
Lemma \ref{prop:odd} and Proposition~\ref{cor:odd} below) that one can 
indeed always find such real $(a,b)$. More to the point, we have reduced the 
study of $Z^*_\R(g)$ from $6$ to $2$ parameters. \dia 
\end{ex} 
\begin{lemma}
\label{prop:odd} 
Suppose $\cA\!=\!\{a_1,\ldots,a_m\}\!\subset\!\Zn$ affinely generates $\Zn$ 
and $a_1\!=\!\bO$. Then there are $i_1,\ldots,i_n\!\in\!\{2,\ldots,m\}$ such 
that $\det[a_{i_1},\ldots,a_{i_n}]$ is odd. \qed 
\end{lemma} 

\begin{dfn} 
\label{dfn:red} 
Suppose $\cA\!=\!\{a_1,\ldots,a_m\}\!\subset\!\Zn$ affinely generates $\Zn$, 
has cardinality $m\!\geq\!n+2$, and $a_1\!=\!\bO$. 
We call any set $C\!=\!\{i_1,\ldots,i_n\}$ 
with $\det[a_{i_1},\ldots,a_{i_n}]$ odd as in Lemma \ref{prop:odd} above, 
an {\bf odd cell} of $\cA$. For any $n\times m$ 
matrix $B$, we then let $B_C$ (resp.\ $B_{C'}$) denote the submatrix of 
$B$ defined by columns of $B$ with index in $C$ (resp.\ $\{2,\ldots,m\}
\setminus C$). For any vectors $v,w\!\in\!(\Cs)^m$, let us denote 
their coordinate-wise product by $v\cdot w\!:=\!(v_1w_1,\ldots,v_m w_m)$.  
Also let $\Gamma$ be the multivalued\footnote{The multiple values arise 
from the presence of rational exponents, and the number of images of a 
point is always bounded above by a constant depending only on $\cA$.} 
function from $(\Cs)^m$ to $(\Cs)^{m-n-1}$ defined by 
$\Gamma(y):=\frac{y_{C'}}{y_1}\cdot \left(\frac{y_C}{y_1}
\right)^{-A^{-1}_C A_{C'}}$. Finally, we define the {\bf reduced} 
$\cA$-discriminant variety, $\bla\!\subset\!\C^{m-n-1}$, to be the closure 
of\\
\mbox{}\hfill$\left\{ \Gamma(u) \; \left| \; u\!:=\!(u_1,
\ldots,u_m)\!\in\!(\Cs)^m,\ \cA u \!=\!\bO, \ \sum^m_{i=1}u_i\!=\!0\right. 
\right\}$, \hfill\mbox{}\\  
and call any connected component of $(\Rs)^{m-n-1}\setminus\bla$ a 
{\bf reduced ($\cA$-)discriminant chamber}. \dia 
\end{dfn} 
\begin{rem}\label{rem:red}
\scalebox{.98}[1]{Since we always implicitly assume that an odd cell has 
been fixed a priori for our}\linebreak 
\scalebox{.9}[1]{reduced $\cA$-discriminant varieties, $\Gamma$ 
in fact restricts to a single-valued function from 
$(\Rs)^m$ to $(\Rs)^{m-n-1}$. \dia} 
\end{rem}
\begin{ex} 
\label{ex:param} 
Continuing Examples \ref{ex:big}, \ref{ex:big2}, and \ref{ex:big3}, let us 
shift our original $\cA$ slightly to instead work with $\cA\!=\!
\text{\scalebox{1}[.7]{$
\begin{bmatrix} 
0 & -6 & -6 & -6 & -3 & -5\\
0 & 3 & 1 & 6 & 0 & 0 \\
0 & 0 & 0 & 1 & 1 & 1 
\end{bmatrix}$}}$. 
(The underlying $\cA$-discriminants are left unchanged thanks to 
affine invariance \cite[Prop.\ 1.4, Pg.\ 272]{gkz94}.) 
As observed in Example \ref{ex:big3}, 
$C\!=\!\{3,4,6\}$ is an odd cell 
for $\cA$. So $A_C$ is then the $3\times 3$ matrix from 
Example \ref{ex:big2}, 
$A_{C'}\!=\!\text{\scalebox{1}[.7]{$\begin{bmatrix}
-6 & -3 \\
 3 &  0 \\
 0 &  1 
\end{bmatrix}$}}$, $A^{-1}_C A_{C'}\!=\!\text{\scalebox{1}[.7]
{$\begin{bmatrix}
33/35 & -12/35 \\
12/35 & 2/35\\
-12/35 & 33/35\end{bmatrix}$}}$,  
and thus $\Gamma(u)\!=\!\left(\frac{u_2u^{12/35}_6}
{u^{2/35}_1u^{33/35}_3u^{12/35}_4},\frac{u_5u^{12/35}_3}
{u^{12/35}_1u^{2/35}_4u^{33/35}_6}\right)$.  
Defining $\ell_1,\ldots,\ell_6$ \mbox{respectively} as the polynomials 
$-2\lambda-2, 35\lambda+11,-33\lambda-9,-12\lambda-4,4,12\lambda$, 
and letting $\Psi(\lambda)\!:=\!(\psi_1(\lambda),\psi_2(\lambda))\!:=\!
\Gamma(\ell_1(\lambda),\ldots,\ell_6(\lambda))$, we thus obtain 
that $\bla$ is the closure of\linebreak 
$\{\Psi(\lambda)\; | \; \lambda\!\in\!\C, 
\ \ell_1(\lambda)\cdots\ell_6(\lambda)\!\neq\!0\}$. \dia
\end{ex} 

Let us now consider how the topology of $Z^*_\R(f)$ changes as 
$f$ ranges through $\cA$-discriminant chambers. First, note that 
it is easy to show that for any $Z\!\subseteq\!\Rsn$ and any 
coordinate reflection $\sigma : \Rn\!\longrightarrow\!\Rn$, $Z$ and $\sigma(Z)$ 
are not only diffeomorphic but diffeotopic. With just a little more 
work, one can then show the following: 
\begin{prop}
\label{cor:odd} 
Suppose $\cA\!=\!\{a_1,\ldots,a_m\}\!\subset\!\Zn$ affinely generates $\Zn$,  
has cardinality $m\!\geq\!n+2$, 
and $a_1\!=\!\bO$. Also let $C$ be any odd cell of $\cA$, 
let $f(x)\!:=\!\sum^m_{i=1}\delta_ix^{a_i}$ with 
$\delta\!:=\!(\delta_1,\ldots,\delta_m)\!\in\!(\Rs)^m$, and let 
$\bar{\delta}\!\in\!(\Rs)^m$ be 
the unique vector with $\bar{\delta}_1\!=\!1$, 
$\bar{\delta}_C\!=\!(1,\ldots,1)$ and 
$\bar{\delta}_{C'}\!=\!\Gamma(\delta)$. Finally, 
let $\bar{f}\!:=\!\sum^m_{i=1}\bar{\delta}_ix^{a_i}$ and 
let $\conv \cA$ denote the convex hull of $\cA$. Then: 
\begin{enumerate} 
\item{ $\Gamma$ induces a surjection from the set of 
connected components of\\ 
\mbox{}\hfill $\Pro^{m-1}_\R\setminus\left(\nabla_\cA\cup
\left\{\left.[y_1:\cdots:y_m]\!\in\!\Pro^{m-1}_\R\; \right| 
\; y_1\cdots y_m\!=\!0 \right\}\right)$
\hfill\mbox{}\\ 
to the set of reduced $\cA$-discriminant chambers. }  
\item{ If, for all facets $Q'$ of $\conv \cA$, we have that 
$\#(\cA\cap Q')\!=\!n$, then $Z^*_\R(f)$ and $Z^*_\R(\bar{f})$ 
are diffeotopic. Furthermore, for any $f_1$ and $f_2$ with 
$\bar{f}_1$ and $\bar{f}_2$ lying in the 
same reduced $\cA$-discriminant chamber, $Z^*_\R(\bar{f}_1)$ and 
$Z^*_\R(\bar{f}_2)$ are diffeotopic. \qed}
\end{enumerate} 
\end{prop}

\noindent
Proposition~\ref{cor:odd} follows easily from a routine application of 
the Smith normal form and the implicit function theorem. In particular, the 
crucial trick is to observe that exponentiation by $A_C$, when $C$ is an odd 
cell, induces an automorphism of orthants of $\Rsn$. 
Our assumption on the intersection of $\cA$ with the facets of $\conv \cA$ 
ensures that any topological change in the  
zero sets of $f$ and $\bar{f}$ (in the underlying real toric variety 
corresponding to $\conv \cA$ \cite{tfulton}) occurs within $\Rsn$. 
\begin{rem}
It is also easily checked that our genericity assumption (on the 
intersection of $\cA$ with the facets of $\conv A$) implies that 
(a) the nullspace of $\cA$ is {\bf not} contained 
in any coordinate hyperplane, and (b) when $\#\cA\!=\!n+3$, one has 
that $\overline{\nabla}_{\cA}\cap\R^2$ is not contained in any 
line. \dia 
\end{rem} 
\begin{dfn} 
\label{dfn:over} 
Following the notation of Definition~\ref{dfn:red}, 
given any $\delta\!\in\!(\Cs)^m$, let $\delta'$ be the unique 
vector with $\delta'_1\!:=\!1$, $\delta'_C\!=\!1$, and $\delta'_{C'}
\!=\!\delta_{C'}$. We then define the {\bf reduced $\cA$-discriminant} 
to be $\overline{\Delta}_\cA(\delta_{C'})\!:=\!\Delta_\cA(\delta')$. \dia 
\end{dfn}  

\subsection{Going Beyond One Polynomial Via the Cayley Trick} 
\label{sub:cayley}
We will need one last construction in order to apply $\cA$-discriminants 
to systems of equations. 
\begin{dfn} 
Let $e_i$ denote the $i^\thth$ standard basis vector of $\R^\infty$. 
Then, for any $\cA_1,\ldots,\cA_k\!\subset\!\Rn$, we call
$\cay(\cA_1,\ldots,\cA_k)\!:=\!(\cA_1\times\{0\}) \cup 
(\cA_2\times\{e_{n+1}\})\cup\cdots \cup
(\cA_k\times\{e_{n+k-1}\})$ 
the {\bf Cayley embedding} of $(\cA_1,\ldots,\cA_k)$. 
We also define the {\bf Newton polytope} of $f$ to be
$\newt(f)\!:=\!\conv(\supp(f))$ and, for any compact set 
$B\!\subset\!\Rn$ and $w\!=\!(w_1,\ldots,w_n)\!\in\!\Rn$, we
define $B^w$ --- the face of $B$ with inner normal $w$ --- to be 
$\{(x_1,\ldots,x_n)\!\in\!B\; | \; x_1 w_1+\cdots+x_n
w_n \text{ is minimized}\}$. \dia
\end{dfn}
\begin{ex}
Taking $(\cA_1,\cA_2)$ to be the support of any $H_{(a,b,3)}$ in the 
Haas family, observe that $\cay(\cA_1,\cA_2)$ is exactly the 
$\cA$ from Examples \ref{ex:big}--\ref{ex:big3}. \dia 
\end{ex}  
\begin{lemma}
\label{lemma:chamber} 
(See, e.g., \cite[Prop.\ 1.7, Ch.\ 9 and Pg.\ 380]{gkz94}.)
Suppose $f_1,f_2\!\in\!\R[x,y]$ are bivariate polynomials with
respective supports $\cA_1$ and $\cA_2$, and let $F\!:=\!(f_1,f_2)$. 
Assume also that for all $w\!\in\!\R^2$ the
pair of faces $(\conv(\cA_1)^w,\conv(\cA_2)^w)$ 
{\bf never} consists of two parallel edges. 
Then --- identifying the coefficient of $x^a$ in $f_i$ with the point
$(a,0)$ or $(a,1)$ of $\cay(\cA_1,\cA_2)$ according as $i$ is $1$ or $2$ ---
the number of isolated real roots of
$F\!:=\!(f_1,f_2)$ in any open quadrant is {\bf constant} on
any discriminant chamber of $\cay(\cA_1,\cA_2)$. 
\end{lemma}

\noindent
{\bf Proof:} Fixing an ordering on $\cA_1$ and $\cA_2$, let $c^{(i)}$ denote 
the coefficient vector of any $f_i$ with support $\cA_i$. 
Also let $\nabla_{(\cA_1,\cA_2)}$ denote the closure in 
$\Pro^{\#\cA_1+\#\cA_2-1}_\C$ of those $\left[c^{(1)}:c^{(2)}\right]$ such 
that the corresponding $F\!:=\!(f_1,f_2)$ has a degenerate root in $(\Cs)^2$. 
Via \cite[Prop.\ 1.7, Ch.\ 9 and Pg.\ 380]{gkz94}, $\nabla_{(\cA_1,\cA_2)}$ 
is exactly $\nabla_{\cay(\cA_1,\cA_2)}$. 

Note also that by our 
assumption on the $(\conv(\cA_1)^w,\conv(\cA_2)^w)$ (using 
$w\!\in\!\{e_1,e_2\}$), 
the only possible isolated root of $F$ on the coordinate cross is $(0,0)$. 
To conclude, we need only to observe that if 
$\left[c^{(i)}\right]\!\in\!\Pro^{\#\cA_i-1}_\C
\setminus\nabla_\cA$ and $c^{(i)}$ has {\bf no} zero coordinates 
for all $i\!\in\!\{1,2\}$,  
then all the roots of $F$ in the toric variety corresponding 
to $\conv(\cA_1)+\conv(\cA_2)$ lie in $(\Cs)^2$. Thus, along any fixed 
path within any fixed $\cay(\cA_1,\cA_2)$-discriminant chamber, the roots of 
$F$ are continuous functions (with bounded range) of the coefficients. \qed 

\section{The Proof of Theorem~\ref{thm:disc} } 
\label{sec:pfdisc} 
The reduced $\cA$-discriminant variety has many interesting properties 
that we will exploit. Before proving Theorem~\ref{thm:disc}, however, 
we will need an important recent bound on the number of real roots of 
certain structured polynomial systems. 
\begin{dfn} 
\label{dfn:shear} 
Suppose $\ell_1,\ldots,\ell_j\!\in\!\R[\lambda_1,\ldots,
\lambda_k]$ are polynomials of degree $\leq\!1$. We then call any system of 
equations of the form 
$S\!:=\!\left(1-\prod^j_{i=1} \ell^{b_{1,i}}_i(\lambda_1,\ldots,\lambda_k),
\ldots,1-\prod^j_{i=1} \ell^{b_{k,i}}_i(\lambda_1,\ldots,\lambda_k)\right)$, 
with $b_{i,i'}\!\in\!\R$ for all $i,i'$, and the vectors 
$(b_{1,1},\ldots,b_{1,j}),\ldots,(b_{k,1},\ldots,b_{k,j})$ linearly  
independent, a {\bf $k\times k$ sheared binomial system with $j$ factors}. 
We also call each $\ell_i$ a {\bf factor} of the system. A sheared binomial system
is referred to as  a
{\bf Gale Dual System} in \cite{bs}.\dia 
\end{dfn} 

\noindent
Note that our definition implies that $j\!\geq\!k$. For $j\!=\!k$, 
it is easy to reduce any $k\times k$ sheared binomial system with $j$ factors 
to a $k\times k$ linear system, simply by multiplying and dividing equations 
(mimicking Gaussian elimination). For $j\!>\!k$, sheared polynomial systems 
become much more complicated. 
\begin{thma}
\label{thm:shear} The number of non-degenerate roots 
$\lambda\!\in\!\R^k$ of any $k\times k$ sheared binomial system with $n+k$ 
factors, and all factors positive, is bounded above by:\\ 
1.\ \cite[Lemma 2]{tri} $n+1$ (and the same bound applies if we also count 
degenerate isolated roots with all factors positive), for $k\!=\!1$,\\ 
2.\ \cite{bs} $(e^2+3) 2^{(k-4)(k+1)/2} n^k$, for all $k\!\geq\!1$. \\
In particular, $e^2+3\!\approx\!10.38905610$. \qed 
\end{thma} 
\medskip

\noindent
We will also need one last 
important result before our main proof.  
Let $\Psi (\lambda) = \Gamma (\ell_1(\lambda), \dots, \ell_m(\lambda))$ be the 
dense parametrization of a reduced $\cA$-discriminant variety associated to an 
odd cell of $\cA$ as in Definition~\ref{dfn:red}. We can in fact consider this 
map to be defined over $\Pro^{1}_\C$ by instead working with the 
homogenizations 
$\ell_i([\lambda_1: \lambda_2] )\!:=\!\alpha_i \lambda_1 + \beta_i \lambda_2$.
Also, recalling Definition \ref{dfn:over} and the fact that 
$\overline{\Delta}_{\cA}$ is an irreducible defining polynomial for $\bla$ 
with integer coefficients, let $Z$ be the finite set of points 
of $\bla$ at which the gradient of $\overline{\Delta}_{\cA}$ vanishes. 
\begin{lemma}
\label{lemma:funk}
With the notation above, if $\Psi([\lambda_1, \lambda_2]) \in (\Rs)^{2} \cap 
\left(\bla \setminus Z\right)$ then $[\lambda_1: \lambda_2]$ can be chosen in 
$\Pro^{1}_\R$. In other words, $\left(\R^2 \cap \bla\right)
\setminus\Psi(\Pro^1_\R)$ is finite.  
\end{lemma}

\noindent 
{\bf Proof of Lemma~\ref{lemma:funk}:} 
The proof is an easy consequence of the fact that $\Psi$ defines a multivalued 
function (univalued from the real points, as we have already noted in 
Remark~\ref{rem:red}) which
is an inverse to the logarithmic Gauss map $G:(\C^*)^2 \cap (\bla \setminus Z) 
\longrightarrow \Pro^1_\C$:
\[G(y):=\left[y_1 \frac{\partial} {\partial{y_1}}\overline{\Delta}_{\cA}(y): 
y_2 \frac{\partial}{\partial{y_2}}\overline{\Delta}_{\cA}(y)\right].\]
This is proved as in~\cite{cd, pt}. Now, since $G$ has rational coefficients,
$G(y)$ has real coordinates for each real point $y\in (\Rs)^{2} \cap 
(\bla \setminus Z)$. \qed 

\medskip 
We are now ready to prove Theorem~\ref{thm:disc}.

\medskip 
\noindent 
{\bf Proof of Theorem~\ref{thm:disc}:}  
Let $\cT_\cA$ denote the toric variety corresponding to the convex 
hull of $\cA$ \cite{tfulton}. Note that by our assumption that
every facet of $\cA$ corresponds to the vertices of a simplex,  
the complex zero set of any $f$ with $\supp(f)\!=\!\cA$ is thus always 
nonsingular at infinity, relative to $\cT_\cA$ (see, e.g., 
\cite[Sec.\ 3.2]{thresh}). 
By Proposition~\ref{cor:odd}, it then suffices to show that our desired 
bound applies to the number of {\bf reduced} $\cA$-discriminant 
chambers. Note also that by Proposition~\ref{cor:odd} and 
Lemma~\ref{lemma:funk}, 
the real part of the reduced $\cA$-discriminant variety --- 
$\R^2\cap\bla$ --- must be the union of a finite set of points and the 
closure of $\{\Psi(\lambda)\; | \; \lambda\!\in\!\R,\ 
\ell_1(\lambda)\cdots\ell_{n+3}(\lambda)\!\neq\!0\}$, 
where $\Psi(\lambda)\!:=\!(\psi_1(\lambda),\psi_2(\lambda))
\!:=\!\left(\prod^{n+3}_{i=1}\ell^{b_{1,i}}_i(\lambda),
\prod^{n+3}_{i=1}\ell^{b_{2,i}}_i(\lambda)\right)$, 
and $\ell_1,\ldots,\ell_{n+3}$ are univariate polynomials in $\lambda$ 
of degree $\leq\!1$ 
Let $\Omega\!\subset\!\R^2$ denote the 
aforementioned closure. Since isolated points do not disconnect connected 
components of the complement of a (locally closed) real algebraic curve, 
it thus suffices to focus on $\Omega$. In particular, the 
connected components of $(\Rs)^2\setminus\Omega$ are (up to 
the deletion of finitely many points) exactly the reduced 
$\cA$-discriminant chambers. Note also, by observing the poles of the 
$\psi_i$, that $\Omega$ is the closure of the union of 
no more than $n+4$ arcs, i.e., homeomorphic images of the open interval 
$(0,1)$.

To count the number of connected components of $(\Rs)^2\setminus\Omega$, we 
will use the classical {\bf critical points method} \cite{chigo}, combined 
with our more recent tools. In particular, let us first bound the number of 
$x$-axis intersections, 
cusps, vertical tangents and vertical asymptotes, and nodes of $\Omega$ in 
$\R^2$. (These 
constitute all possible critical points of the orthogonal projection mapping 
$\Omega$ to the first coordinate axis.) Let us call the 
numbers of these respective objects $M_0(\cA)$, $M_1(\cA)$, $M_2(\cA)$, and 
$M_3(\cA)$, and proceed with bounding their number from above. 

\medskip
\noindent 
\fbox{$x$-axis Intersections:} Clearly, $\Omega$ intersects the 
$x$-axis iff $\vp_2(\lambda)\!=\!0$, and the latter occurs iff a monomial in 
the $\ell_i$ vanishes at some $\lambda\!\in\!\C\cup
\{\pm\infty\}$. Also, by Theorem~\ref{thm:dfs} and Definition~\ref{dfn:red}, 
the degrees of the numerator and denominator of $\vp_2$ are equal.  
So there are clearly no more than $n+2$ solutions to $\vp_2(\lambda)\!=\!0$, 
and thus $M_0(\cA)\!\leq\!n+2$. 

\medskip
\noindent 
\fbox{Cusps (and certain isolated real points):} 
To count cusps, it suffices to bound from above the number of 
complex $\lambda$ such that $\frac{\partial \vp_1}{\partial \lambda}
\!=\!\frac{\partial \vp_2}{\partial \lambda}\!=\!0$. (Note also that 
for such a $\lambda$, $\Psi(\lambda)$ could also be an isolated real 
point of the real part of $\overline{\nabla}_\cA$ off of $\Omega$.) 
In particular, via the product rule, and by \linebreak dividing out by 
suitable monomials in $\ell_1(\lambda),\ldots,\ell_{n+3}(\lambda)$, 
the preceding equation reduces to a univariate system of the 
form $\sum^{n+3}_{i=1}\frac{b_{1,i}\ell'_i(\lambda)}{\ell_i(\lambda)}=
\sum^{n+3}_{i=1}\frac{b_{2,i}\ell'_i(\lambda)}{\ell_i(\lambda)}=0$.  
We can then multiply through by $\ell_1(\lambda)\cdots\ell_{n+3}(\lambda)$ 
to obtain a univariate polynomial of degree $n+2$. Furthermore, it 
is easily checked that the maximal number of distinct cusps as 
$\lambda\longrightarrow\pm\infty$ is one, 
and thus $M_1(\cA)\!\leq\!n+3$.  

\medskip
\noindent 
\fbox{Vertical Tangents (and vertical asymptotes):} Here, we proceed 
essentially the same 
as for cusps, but with only one derivative. However, 
$\lambda\!\in\!\{\pm\infty\}$ can possibly yield two distinct vertical 
tangents. So $M_2(\cA)\!\leq\!n+4$. Note also that this count 
includes all cusps.  

\medskip
\noindent 
\fbox{Nodes:} 
Here, we need to bound 
the number of $2$-sets $\{\lambda,\lambda'\}$ (so $\lambda\!\neq\!\lambda'$) 
with $\Psi(\lambda) \!=\!\Psi(\lambda')$ and 
$\lambda,\lambda'\!\in\!\R\cup\{-\infty,+\infty\}$. 
Counting these real pairs then reduces to counting the number 
of real solutions of a $2\times 2$ sheared binomial system --- with 
$\leq\!2(n+3)\!=\!2n+6$ factors --- where no factor is zero. 
Theorem~\ref{thm:shear} 
counts such solution satisfying certain sign condition for 
the $\ell_i$, so let us carefully count the number of (nonzero)  
sign combinations possible for the vectors 
$(\ell_1(\lambda),\ldots,\ell_{n+3}(\lambda))$ and  
$(\ell_1(\lambda'),\ldots,\ell_{n+3}(\lambda'))$: 
Clearly, each such a vector admits at most $(n+4)$ possible (nonzero) sign 
combinations, since the sign of any $\ell_i$ is constant to the right (or 
to the left) of its unique real root, and there are no more than $n+3$ real 
roots for our $\ell_i$. Thus, there are at most $(n+4)^2$ possibilities 
for the (nonzero) sign vector of $(\ell_1(\lambda),\ldots,\ell_{n+3}(\lambda),
\ell_1(\lambda),\ldots,\ell_{n+3}(\lambda))$. 
So, combining with Theorem~\ref{thm:shear}, 
and noting that there are $\leq\!(2n+4)+2$ factors, we thus 
clearly obtain no more than 
$(n+4)^2 \cdot (e^2+3)\cdot 2^{-2}\cdot (2n+4)^2$ pairs $(\lambda,\lambda')
\!\in\!\R^2$ with $\Psi(\lambda)\!=\!\Psi(\lambda')$ and 
$\lambda\!\neq\!\lambda'$. Note also that there are infinitely 
many solutions of $\Psi(\lambda)\!=\!\Psi(\lambda')$ of the form 
$\lambda\!=\!\lambda'$, but these are non-isolated and thus {\bf not} counted 
by Theorem~\ref{thm:shear}. So there are no more than 
$\frac{e^2+3}{2}(n+2)^2(n+4)^2$ nodes arising from $2$-sets in 
$\R^2$, since our underlying sheared system is symmetric. 

Now, should a $\lambda\!\in\!\R$ 
yield $\Psi(\lambda)\!=\!\Psi(\pm \infty)$ as a node (for some 
fixed choice of sign), then $\vp_2(\lambda')$ must clearly 
have a well-defined nonzero limit as $\lambda'\longrightarrow \pm\infty$,  
since we have already counted $x$-axis intersections and vertical tangents. 
We are thus reduced to counting the number of real roots of a univariate 
sheared binomial, with no factor zero. By Theorem~\ref{thm:shear} again 
(with $\leq\!(n+2)+1$ factors), 
and recalling our last observation on the sign vector of $(\ell_1,\ldots,
\ell_{n+3})$, we then directly obtain no more than $(n+4)(n+2)$ 
nodes arising from $(\lambda,\lambda')\!\in\!(\R\cup\{\pm\infty\})^2\setminus
\R^2$. (It is also easily checked that should $(\pm\infty,\mp\infty)$ yield a 
node, then there is another pair $(\lambda,\lambda')\!\in\!\R^2$ giving 
the same node.) 

In summary, we thus obtain no more than 
$\frac{e^2+3}{2}(n+2)^2(n+4)^2+(n+2)(n+4)$ nodes, and thus   
$M_3(\cA)\!\leq\!(n+2)(n+4)\left[1+\frac{e^2+3}{2}\cdot (n+2)(n+4)\right]$.  
  
\medskip
\noindent 
\fbox{Back to Reduced Chambers...} 
To count the number of connected components of 
$(\Rs)^2\setminus\Omega$, let us now introduce vertical lines 
$L_1,\ldots,L_N$ exactly at the locations of the $y$-axis, the $x$-axis 
intersections, the cusps, the vertical tangents, and the nodes 
of $\Omega$. Clearly, any connected component of\\ 
\mbox{}\hfill$T\!:=\!(\Rs)^2\setminus(\Omega\cup L_1\cup \cdots\cup L_N)$
\hfill\mbox{}\\ 
is contained in a unique connected component of $(\Rs)^2\setminus\Omega$. So 
it suffices to count the connected components of $T$. To do the latter, 
observe that $N\!\leq\!1+M_0(\cA)+M_2(\cA)+M_3(\cA)$ (since 
cusps were already counted among vertical tangents via our technique 
above), and that our lines $\{L_i\}$ thus divide $(\Rs)^2$ into no more than 
$1+N\!=\!1+1+(n+2)+(n+4)+(n+2)(n+4)\left[1+\frac{e^2+3}{2}\cdot (n+2)(n+4)
\right]$\linebreak 
$=16+8n+n^2+\frac{e^2+3}{2}\cdot (n+2)^2(n+4)^2$ vertical strips. 

Now note that within the interior of each strip, $\Omega$ 
is smooth, does {\bf not} intersect the $x$-axis, and has no vertical 
tangents. So to count components of $T$ within 
any particular vertical strip, we need only bound from above  the number of 
non-degenerate intersections of $\Omega\cup\{x_2\!=\!0\}$ with a vertical 
line distinct from $L_1,\ldots,L_N$. This clearly reduces to counting the 
number of real roots of a binomial 
in $n+3$ univariate linear forms. Via our earlier sign condition count, 
and by Theorem~\ref{thm:shear} once again (with $\leq\!(n+2)+1$ factors), the 
desired upper bound is then $1+(n+4)(n+2)$. 
Thus, each of our vertical strips contains no more than 
$2+(n+2)(n+4)$ connected components of $T$. Taking into account the number of 
vertical strips, we thus finally arrive at an upper bound of\\  
\mbox{}\hfill$\left(10+6n+n^2\right)\left(16+8n+n^2+\frac{e^2+3}{2}\cdot 
(n+2)^2(n+4)^2\right)$\hfill \mbox{}\\ for the number 
of connected components of $T$, our bound is proved. 

To conclude, note that Assertion (b) follows immediately 
from our bound, and Assertion (a) follows from merely comparing 
coefficients in the underlying polynomials. \qed 
\begin{rem} 
The diffeotopies we have used above are thus obtained 
by (i) following a path in a reduced $\cA$-discriminant chamber, and 
(ii) performing a coordinate reflection. One could certainly allow a 
broader class of diffeotopies, and thus (potentially) greatly reduce 
the bound we have just proved. Diffeotopies obtained solely via (i) are 
analogues of what are sometimes called {\bf rigid} diffeotopies in other 
settings (see, e.g., \cite{orevkov}). \dia 
\end{rem}  
\begin{rem} 
Isolated points {\bf can} in fact occur in the real part of a 
reduced $\cA$-discriminant\linebreak variety. For instance, taking 
$\cA\!=\!\text{\scalebox{1}[.7]{$\begin{bmatrix}
0 & -6 & -6 & -6 & -3 & -5\\
0 & 3 & 1 & 6 & 0 & 0 \\
0 & 0 & 0 & 1 & 1 & 1 \end{bmatrix}$}}$ (as in Example \ref{ex:param}), 
we will see in the next section (and the Appendix) that 
$\R^2\cap\overline{\nabla}_{\cA}$ is the disjoint union of a connected 
finite union of smooth arcs and exactly $3$ isolated points 
(located respectively in the $+-$, $--$, and $-+$ quadrants). This 
in turn implies that the {\bf real part} of the (non-reduced) discriminant 
variety $\nabla_{\cA}$ has connected\linebreak
\scalebox{.95}[1]{components of codimension at least $2$, 
even though $\nabla_{\cA}$ is itself codimension $1$ and irreducible over 
$\C$. \dia} 
\end{rem} 

\section{The Proof of Theorem~\ref{thm:new} } 
\label{sec:new} 
Before going into our main proof, let us first review an important 
criterion for an ``approximate'' root to converge quickly under 
Newton iteration to a true root of a polynomial system. 
\begin{dfn}
\label{dfn:ops} 
\cite{smale,bcss}
Given any analytic function $F : \Cn \longrightarrow \Cn$, 
we let $F'$ denote its Jacobian matrix, and define the 
{\bf Newton endomorphism}, $N_F :  \C^n \longrightarrow \C^n$ to 
be the function $N_F(z)\!:=\!z-F'(z)^{-1}F(z)$. Also, given 
any $z_0\!\in\!\Cn$, we define the sequence of 
{\bf Newton iterates of $z_0$ (under $F$)} to be 
$(z_n)_{n\in\N\cup\{0\}}$ where $z_{n+1}\!:=\!N_F(z_n)$ for all $n\!\geq\!0$. 
Finally, given any multi-linear operator $\cL : (\C^N)^k\longrightarrow 
\C^N$ and $v\!=\!(v_1\ldots,v_N)\!\in\!\C^N$, we let $|v|\!:=\!
\sqrt{|v_1|^2+\cdots +|v_N|^2}$ denote the usual Hermitian norm and let 
$|\cL|$ be the multi-linear operator norm 
$\max\limits_{(v^1,\ldots,v^k)\in (\C^N\setminus \bO)^k 
}\frac{\left|\cL(v^1,\ldots,v^k)\right|}{|v^1|\cdots
|v^k|}$. \dia 
\end{dfn} 
\begin{dfn}
\label{dfn:approx} 
\cite{smale,bcss}  
Following the notation of Definition~\ref{dfn:ops}, 
we define the invariants $\beta(F,z)\!:=\!\left|z-N_F(z)\right|\!=\!
\left|F'(z)^{-1}F(z)\right|$, $\gamma(F,z)\!:=\!\sup_{k\geq 2}
\left|\frac{1}{k!}F'(z)^{-1}F^{(k)}(z)\right|^{1/(k-1)}$, and 
$\alpha(F,z)\!:=\!\beta(F,z)\gamma(F,z)$. Also, 
let us call a point $z_0\!\in\!\Cn$ an {\bf approximate root 
of $F$} iff the Newton iterates of 
$z_0$ under $F$ satisfy $|\zeta-z_n|\!\leq\!\left(\frac{1}{2}\right)^
{2^{n-1}}|\zeta-z_0|$ for all $n\!\geq\!1$, for some 
true (and non-degenerate) root $\zeta\!\in\!\Cn$ of $F$. \dia 
\end{dfn}
\begin{rem}
Note that $F'(z)^{-1}F^{(k)}$ is a symmetric $k$-linear
operator from $(\Cn)^k$ to $\Cn$, and that the underlying coefficients can be
identified with $k^\thth$ order partial derivatives of the $f_i$
constituting $F$. In particular, observe that when $F$ is a
polynomial system, the supremum in the definition of\linebreak 
\scalebox{.98}[1]{$\gamma(F,z_0)$ is in 
fact a maximum over a finite set with cardinality depending on the degrees 
of the $f_i$. \dia} 
\end{rem}

Note that approximate roots (as defined above) give an efficient, rigorous,  
and numerically feasible way to encode true roots: For instance, 
instead of specifying $n$ (likely huge) minimal polynomials for an 
algebraic point $(\zeta_1,\ldots,\zeta_n)$, and $n$ corresponding isolating 
intervals, we can instead simply give an $n$-tuple  
$(z^1_0,\ldots,z^n_0)$ that is an approximate root. One can then extract 
arbitrarily high accuracy through a small number of Newton iteration, thanks 
to Definition~\ref{dfn:approx}.  

There has been much important work on proving useful 
complexity bounds for this approach (see, e.g., \cite{pardo}). In the 
space here, we merely point to \cite{bcss} as an excellent beginning 
reference.  In particular, the $\alpha$ invariant gives a sufficient 
criterion to guarantee that a given point, and any point sufficiently near,  
is an approximate root.  
\begin{thma} 
\label{thm:alpha} 
\cite[Ch.\ 8]{bcss} 
Following the notation above, suppose $z_0\!\in\!\Cn$ satisfies 
$\alpha(F,z_0)\!<\!0.03$. Then $z_0$ --- and any point within distance 
$\frac{.05}{\gamma(F,z_0)}$ --- is an approximate root of $F$. 
Furthermore, the unique root $\zeta$ of $F$ to which the Newton iterates of 
$z_0$ converge satisfies $|z_0-\zeta|\!\leq\!2\beta(F,z_0)$. \qed  
\end{thma} 

\noindent 
The preceding result is sometimes referred to as a {\bf (robust) one 
point estimate}, and considerably strengthens earlier seminal results of 
Kantorovich, from the 1960s,  
which relied on invariants defined over entire regions. We are now ready to 
prove Theorem~\ref{thm:new}. 

\medskip 
\noindent
{\bf Proof of Theorem~\ref{thm:new}:}  Let $\cH(3)$ denote the 
Cayley embedding of the support of any $H_{(a,b,3)}$ in the 
Haas family. Proposition~\ref{cor:odd} and Lemma~\ref{lemma:chamber}  
then tell us that we can find $E_3$ by using the critical points 
method \cite{chigo}, just as in the proof of Theorem~\ref{thm:disc}. In particular, 
for the setting at hand, this reduces to the following four 
steps: 
\begin{enumerate}
\item{computing $R(a,b)\!:=\!\overline{\Delta}_{\cH(3)}(-a,-b)$,  
i.e., $\Delta_{\cH(3)}(1,a,-1,1,b,-1)$ up to sign}    
\item{computing the partial derivatives of $R$ up to order $2$ 
(needed later to check vertical tangents and convexity of the underlying arcs)}
\item{isolating those real $a_0$ such that there is a real $b_0$ with 
$R(a_0,b_0)\!=\!\left.\frac{\partial R}{\partial b}\right|_{(a_0,b_0)}\!=\!0$ }
\item{computing the number of roots of $H_{(a,b,3)}$ in each quadrant, at 
representative choices of $(a,b)$, picking at least one representative 
pair from each reduced discriminant chamber}
\end{enumerate}   
These computations are routine, albeit hours-long, via {\tt Maple}. 
So we now summarize the crucial details. 

\noindent
\fbox{Steps (1) and (2):} 
We calculate $R$ by first computing the resultant $p$ (resp.\ $q$) of 
$h_1$ (resp.\ $h_2$) and the determinant of the Jacobian of $(h_1,h_2)$, 
with respect to the variable $y$, where $(h_1,h_2)\!:=\!H_{(a,b,3)}$. 
We then compute the resultant, $\widetilde{R}\!\in\!\Z[a,b]$, of $p$ and $q$ 
with respect to the variable $x$. $\widetilde{R}$ is a multiple of our 
desired $R$, and we can isolate $R$ by factoring $\widetilde{R}$ and picking 
out the factor of the correct degree. In particular, a quick volume   
calculation (see, e.g., \cite[Thm.\ 4.2.4]{why}) shows that the degree of $R$ 
is bounded above by $236$, and there is only one such factor of 
$\widetilde{R}$. 

It then turns out that $R$ is symmetric, with total degree $90$, degree $47$ 
with respect to each variable, and exactly $58$ monomial terms (each 
having a coefficient with between $43$ and $56$ decimal digits). 
In particular, up to sign, $R(a,b)$ is the polynomial from Example 
\ref{ex:big}, and 
the full monomial term expansion of $R(a,b)$ can be downloaded from  
\hfill{\tt 
http://www.math.tamu.edu/\~{}rojas/haas3disc} . 
\begin{rem} 
\label{rem:comp} 
It is also interesting 
to observe that $R(a,b)$ is of the form $r\!\left(a^6b^{-1},b^6a^{-1}\right)$, 
where $r$ is a polynomial with total degree $18$ and degree $9$ with respect 
to each variable. \dia 
\end{rem} 

\noindent
\fbox{Step (3):} Computing those real $a_0$ where both $R(a_0,b_0)$ and 
$\left.\frac{\partial R}{\partial b}\right|_{(a_0,b_0)}$ vanish 
can be done by computing the resultant of 
$\left(R,\frac{\partial R}{\partial b}\right)$, to eliminate the variable  
$b$. The resulting eliminant $X(a)$ is large, occupying a text file 
over 3Mb in size, but has just $4$ non-monomial factors  
--- $A^*$, $B^*$, $C^*$, and $D^*$, of respective degrees $1260$, $70$, $35$, 
and $5$ --- that have real roots. These factors are respectively 
of the form $A(a^{35})$, $B(a^{35})$, $C(a^5)$, and $D(a^5)$, 
for certain polynomials $A$, $B$, $C$, and $D$ of significantly lower degree: 
$D$ is exactly the minimal polynomial for the rational number 
$\frac{16807}{2916}$, $C$ is the degree $7$ minimal polynomial 
for $\alpha$, and $A$ is the minimal polynomial for the 
real algebraic numbers $\beta$ and $\gamma$ 
(as mentioned in the statement of Theorem~\ref{thm:prob}).  
$B(a)$ is the following quadratic polynomial:\\
\scalebox{.394}[1]
{$682167754844939988089012017644672936127567587235626923090509824
a^{2}
+264684215537625374546748588284030550129153533531664375000000000000000a+
25668846576874695337875882633144450853905456820418476127088069915771484375$.}

The critical values are thus the real roots of $A^*$, $B^*$, $C^*$, $D^*$, and 
there are exactly $16$ of them: $2$ of which are negative and $14$ of which are 
positive. With a bit more work, it is then easily verified that the $-+$  
and $+-$ quadrants of $\R^2$ each contain exactly $2$ reduced 
chambers (each unbounded), and the entire $--$ quadrant is itself 
a reduced chamber. Thus, all the action occurs in the $++$ quadrant, 
and we will need to consider a total of $16+2+1$ vertical strips in our 
application of the critical points method to our real discriminant 
curve $\Omega\!:=\!\R^2\cap\overline{\nabla}_{\cH(3)}$ (two extra vertical 
lines coming 
from the $b$-axis and the unique $a$-axis intersection of $\Omega$).  
The exact location of our vertical strips is detailed in the Appendix. 

There turn out to be exactly $15$ reduced chambers in 
$\R^2_+$: exactly $5$ that are unbounded 
and $10$ that are bounded. An attractive illustration of the unbounded 
reduced chambers can be obtained by computing the logs of the absolute 
values of the coordinates of the zero set of $R$ in $(\Cs)^2$, i.e., the {\bf 
Archimedean amoeba} of $R$ \cite[Cor.\ 1.8]{gkz94}. In particular, the $5$ 
unbounded reduced chambers (and the sole bounded chamber adjacent the origin, 
stretched to {\bf un}boundedness --- by the logarithm --- in the illustration 
below), each contain exactly one of the white convex regions 
below. The remaining bounded chambers thus have images lying inside 
the shaded amoeba, and their boundaries are indicated by the 
darker curve: the image of $\R^2\cap\overline{\nabla}_{\cH(3)}$ under the
log of absolute value map. 

\vspace{.7cm}
\noindent
\scalebox{1}[1]{
\begin{picture}(200,120)(-75,25)
\put(29,-64){\epsfig{file=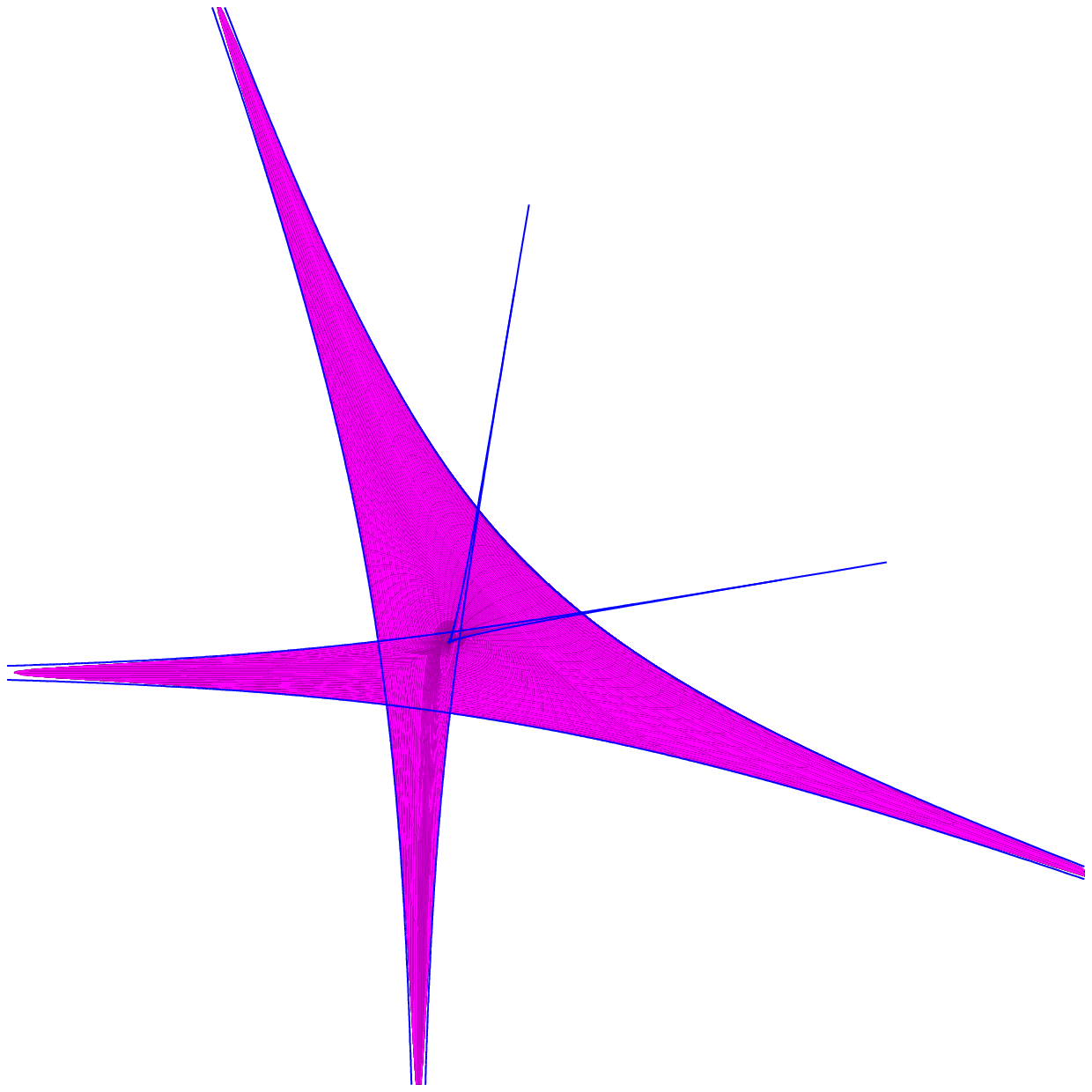,height=3.3in}}
\put(59,-32){\epsfig{file=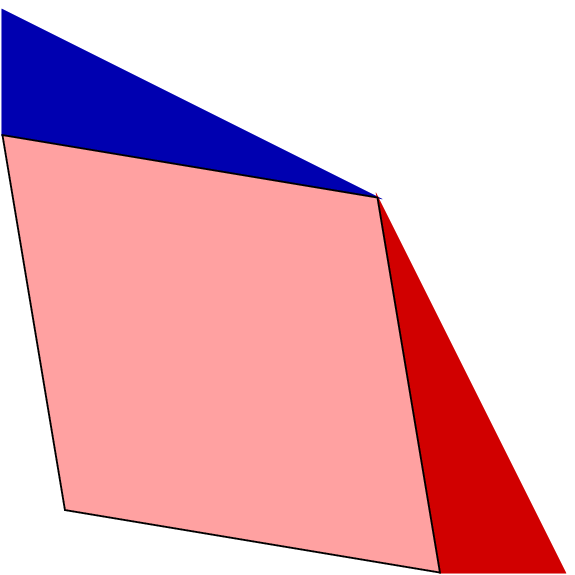,height=.4in}}
\put(164,73){\epsfig{file=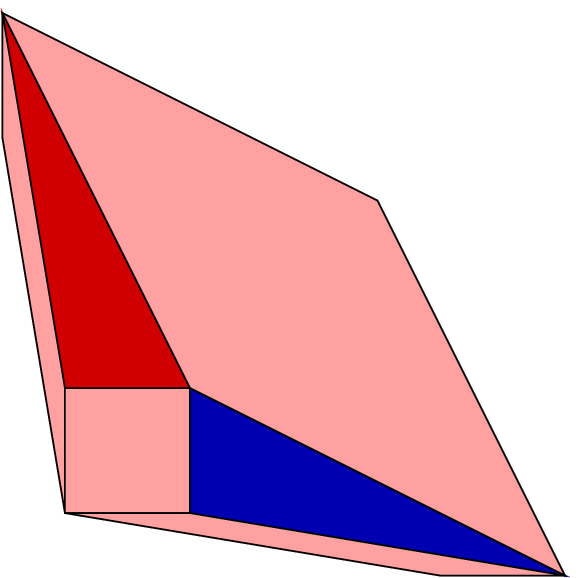,height=.4in}}
\put(106,103){\epsfig{file=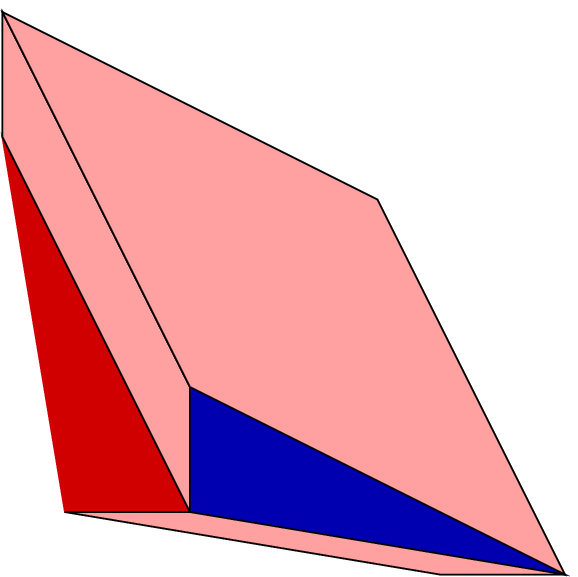,height=.4in}}
\put(196,13){\epsfig{file=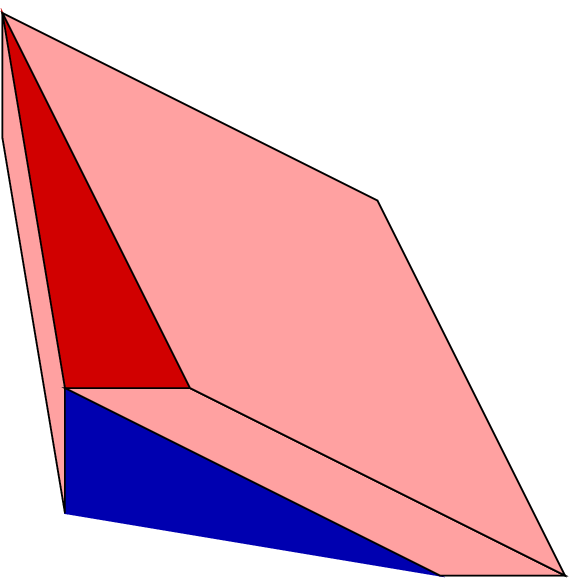,height=.4in}}
\put(61,58){\epsfig{file=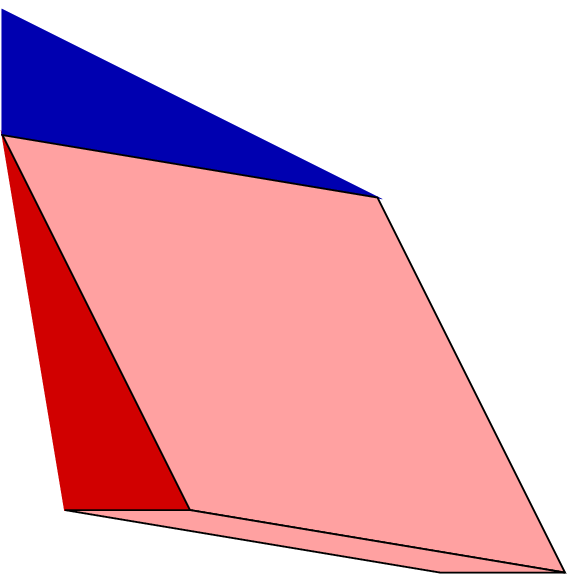,height=.4in}}
\put(151,-38){\epsfig{file=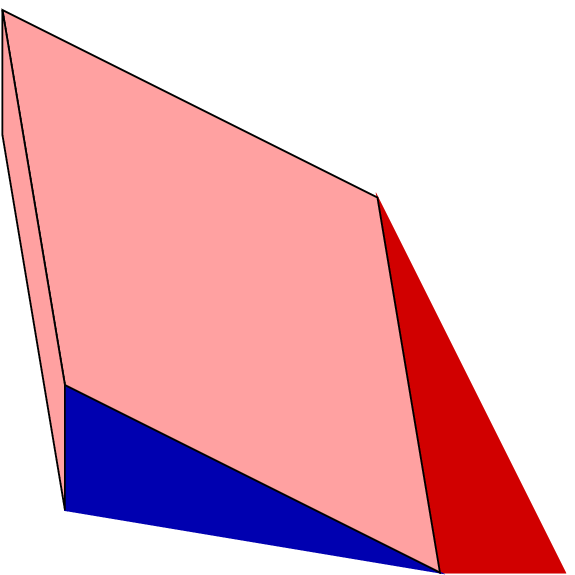,height=.4in}}
\end{picture}}

\vspace{3.1cm}
\begin{rem} 
\label{rem:amoeba} 
The Amoeba Theorem of Gelfand, Kapranov, and Zelevinsky 
\cite[Cor.\ 1.8]{gkz94} tells us 
that the white convex regions above correspond exactly to the vertices 
of $\newt(R)$. The fact that all the facets of the convex 
hull of $\cH(3)$ are triangles 
implies that $\newt(R)$ is actually the image of the {\bf secondary 
polytope} of $\cH(3)$ under an injective linear map \cite[Thm.\ 1.7 of 
pg.\ 221 and Thm.\ 1.4 of pg.\ 302]{gkz94}. Concretely, 
this means that the white convex regions also correspond exactly 
to the triangulations of $\cH(3)$, which we have drawn above as 
well.\footnote{Truthfully, we drew the 
{\bf mixed subdivisions} of the pair $(\supp(h_1),\supp(h_2))$, where 
$H_{(a,b,3)}\!=\!(h_1,h_2)$. The latter diagrams can be interpreted 
as projections of the triangulations of $\cH(3)$.} \dia 
\end{rem} 

\noindent
\fbox{Step (4):} 
Now let $T^*$ be the complement 
of $\Omega\cup L$ in $(\Rs)^2$, where 
$L$ is the set of vertical lines located at the real roots of the 
eliminant $X(a)$ computed in Step (3). Clearly then, 
each connected component of $T^*$ is contained in a unique reduced 
$\cH(3)$-discriminant chamber.  

Returning to the positive quadrant, exactly one of the reduced chambers 
there ($E_3$) possesses $H_{(a,b,3)}$ with $5$ 
positive roots: all other chambers result in $4$ or fewer positive roots. 
Verifying this reduces to solving a representative $H_{(a,b,3)}$ for each 
connected component of $T^*$: there are $125$ such systems.  
In particular, after this sampling of representative points $(a,b)$, we obtain 
that $E_3$ is contained in the union of three adjacent vertical strips, with 
extreme end-points located approximately (to $10$ decimal places) at 
$\{1.4176759490,1.4195167977\}$. 
We also obtain that within each strip, $E_3$ is the region between two smooth 
curves.  This is illustrated below. 

\noindent 
\begin{picture}(200,60)(-80,32)
\put(60,-93){\epsfig{file=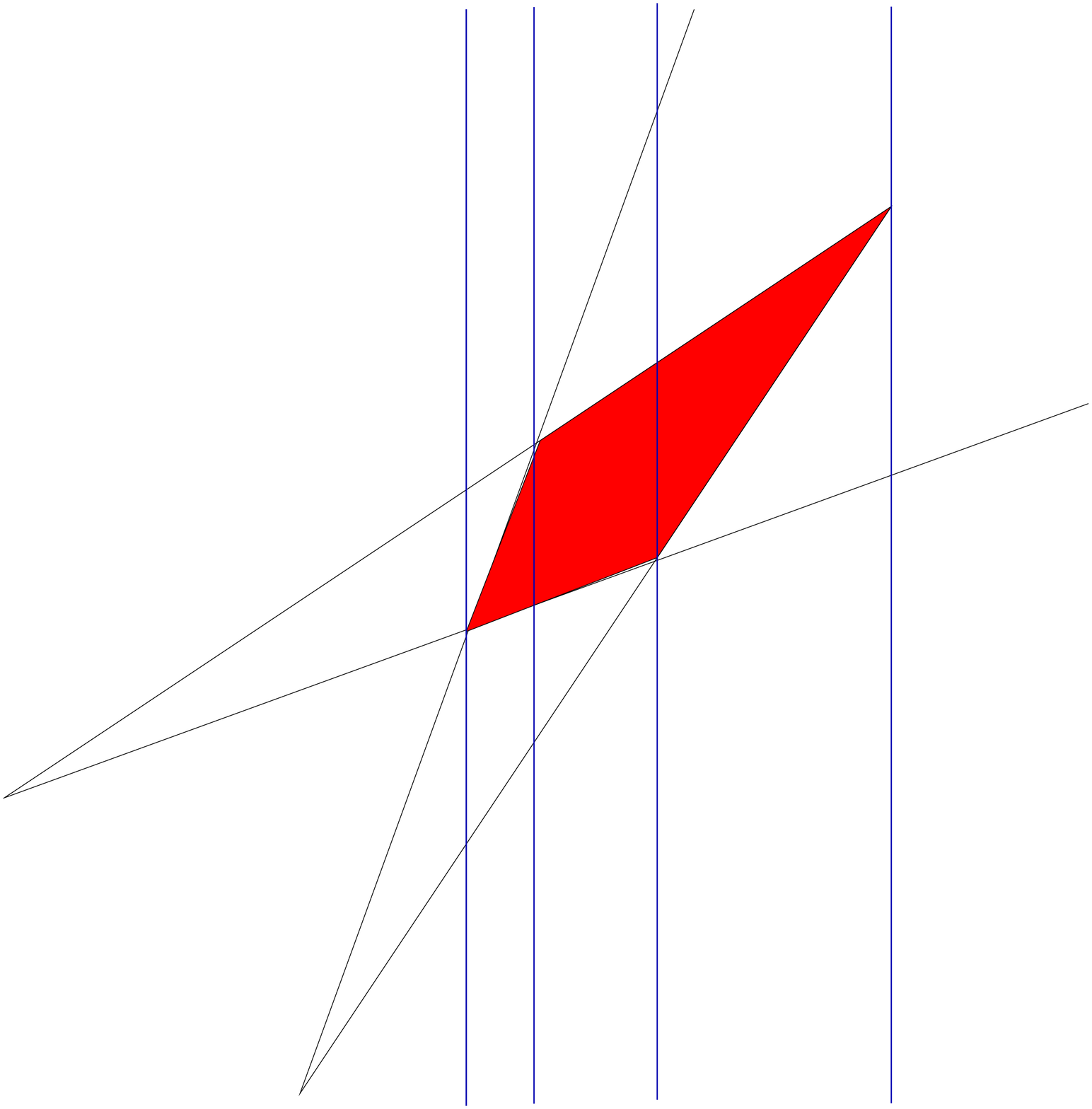,height=2.5in}}
\put(152,6){{\footnotesize $E_3$}} 
\put(-80,-104){{\sc Figure 4.7:} {\em The $(a,b)$ yielding $H_{(a,b,3)}$ with 
$5$ isolated roots in the positive quadrant lie in the}} 
\put(-15,-116){{\em shaded region $E_3$. The vertices 
--- intersected by the four vertical lines --- 
are a}}
\put(-15,-128){{\em subset of the critical points of the projection mapping 
$\overline{\nabla}_{\cH(3)}$ onto the $a$-axis.} } 
\end{picture}

\vspace{5.85cm}
\noindent 
The boundary of $E_3$ thus appears to be 
a subset of the union of $4$ convex arcs. More rigorously, 
it suffices to show that each of the $7$ smooth sub-arcs above, obtained 
from the parametric formula of Example \ref{ex:param}, is convex away 
from its cusps and nodes. This follows easily from computing the partial 
derivatives of the logarithms of our parametric formula and checking 
signs. So $E_3$ is indeed non-empty, star-convex, and in fact has 
$4$ vertices determined by the polynomials mentioned in 
the statement of Theorem~\ref{thm:prob}.  

To conclude, we need only verify that $(44/31,44/31)\!\in\!E_3$. Instead of 
doing this via symbolic algebra, let us 
instead employ Smale's Alpha Theory, as summarized earlier: 
Clearly, we need only check partial 
derivatives of the $h_i$ up to order $6$, and a quick 
computation reveals that any $z\!=\!(z_1,z_2)$ chosen from one of the 
{\bf five} following points satisfies  
$\alpha(H_{(44/31,44/31,3)},z)\!<\!0.03$:\\
\mbox{}\hfill$(\ 0.584513273807\ , \ 0.818672114695\ )$ 
, \ $(\ 0.721441819886\ ,\ 0.757201442567\ )$, \hfill\mbox{}\\ 
\scalebox{.85}[1]{$(0.740238978217,0.740238978217)$}, \ 
\scalebox{.85}[1]{$(0.757201442567,0.721441819886)$}, \ 
\scalebox{.85}[1]{$(0.818672114695,0.584513273807)$}.\linebreak  
Thus, each of these points is an approximate root of $H_{(44/31,44/31,3)}$. 
Also, since the roots of $H_{(44/31,44/31,3)}$ are clearly 
symmetric about the line $\{x\!=\!y\}$, we need only compute 
the $\alpha$-invariant $3$ times. 
So, thanks to Theorem~\ref{thm:alpha} (and Proposition~\ref{prop:alpha} 
of the Appendix), we are done. \qed 
\addtocounter{dfn}{1} 
\begin{rem} 
\label{rem:grobner} 
\scalebox{.94}[1]{Alternatively, we could have simply used any Gr\"obner basis 
solver to get a rational}\linebreak 
\scalebox{.97}[1]{univariate reduction for 
$H_{(44/31,44/31,3)}$. One could then use Sturm-Habicht sequences 
[Stu35, Hab48,}\linebreak
Roy96, LM01] to find (certifiably 
correct) isolating intervals for the real roots. 
However, this naive Gr\"obnerian approach becomes infeasible for higher 
degree examples (cf.\ Note Added in Proof). \dia 
\end{rem}  

\section{The Proof of Theorem~\ref{thm:prob} } 
\label{sec:prob} 
Here, we need only continue the development of the proof of 
Theorem~\ref{thm:new} one step 
further: Since we already observed and proved 
the structure of the boundary of $E_3$ in Section \ref{sec:new},  
we need only verify that $E_1$ and $E_2$ are empty, and 
make two estimates concerning the size of $E_3$. 
The emptiness of $E_2$ follows from essentially the same techniques 
as we used for $E_3$, but the resulting computations (which we omit) are much 
simpler. The emptiness of $E_1$ follows directly from B\'ezout's Theorem: 
$H_{(a,b,1)}$ --- being a pair of bivariate quadratic polynomials --- has no 
more than $4$ non-degenerate isolated {\bf complex} roots. Also, since we 
already know from 
Proposition~\ref{prop:5roots} that $E_3$ is open, and our proof of 
Theorem~\ref{thm:new} already showed $E_3$ to be 
non-empty, we clearly have $\mathrm{Area}(E_3)\!>\!0$.  

To conclude, since the boundary curves of $E_3$ are concave, 
the area of $E_3$ is clearly bounded from above by the area of 
the convex hull of its vertices. {\tt Maple} easily yields the 
estimate stated in our theorem. As for the estimate on the 
probability, we need only observe that the probability is in 
turn bounded above by the aforementioned area, times the 
value of the probability density function at the 
lower left vertex of $E_3$. (The lower left vertex is clearly the point of 
$E_3$ maximizing any function of the form $\alpha e^{\frac{-(x^2+y^2)}
{\beta}}$, when $\alpha,\beta\!>\!0$.) Via {\tt Maple} once again, we are 
done. \qed 

\medskip
Note that our proof here in essence focussed on one reduced 
$\cA$-discriminant chamber, for $\cA\!=\!\text{\scalebox{.7}[.7]
{$\begin{bmatrix}
0 & -6 & -6 & -6 & -3 & -5\\
0 & 3 & 1 & 6 & 0 & 0 \\
0 & 0 & 0 & 1 & 1 & 1
\end{bmatrix}$}}$, while the proof in the last section described a 
collection of sets refining the true reduced $\cA$-discriminant chambers. 
With the calculations we have already done, we can in fact describe {\bf just} 
the relevant $\cA$-discriminant chambers, and give the number of 
real roots for {\bf all} non-degenerate $H_{(a,b,3)}$. In particular, 
the intersection of $\R^2\cap\overline{\nabla}_\cA$ with any sufficiently 
large square is, {\bf up to diffeotopy}, the curve (possessing 
$3$ isolated points) drawn below:\\
\begin{picture}(200,190)(0,25)
\put(0,-5){\epsfig{file=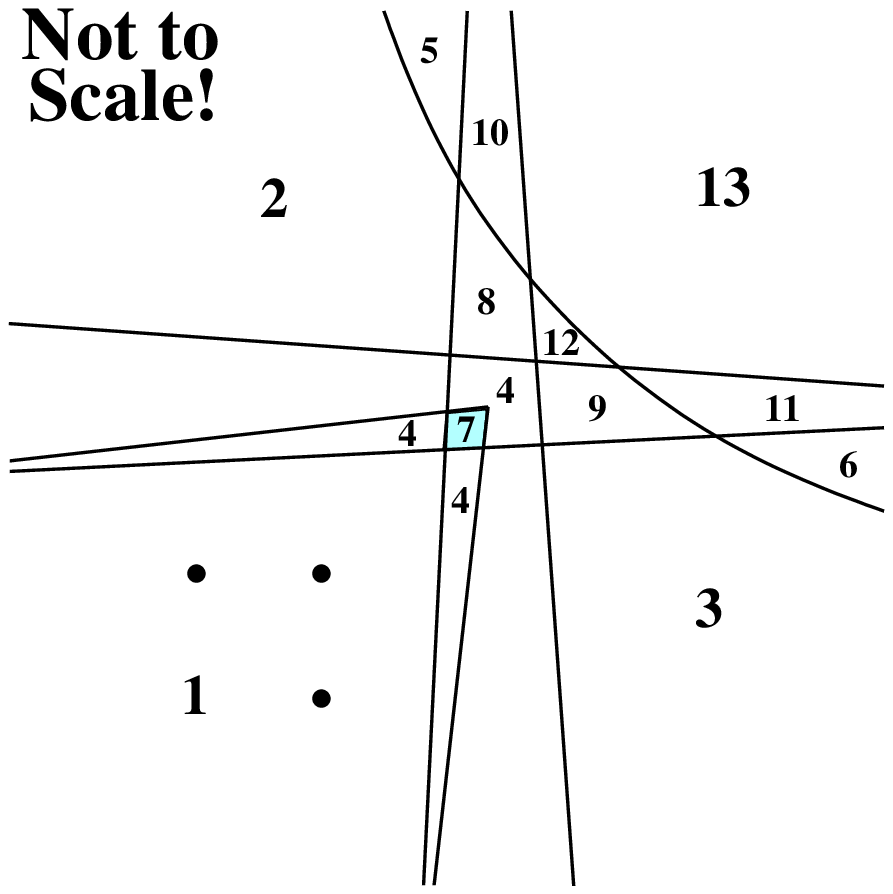,height=2.8in}}
\put(240,110){\scalebox{.9}[.9]{\begin{tabular}{c||c|c|c|c|c}
region \# & $++$ & $-+$ & $--$ & $+-$ & rep.\ $(a,b)$ \\ \hline
$1$ &  $1$ & $0$ & $0$ & $0$ & $(-1,-1)$ \\ \hline
$2$ & $1$ & $0$ & $0$ & $2$ & $(-1,5)$ \\ \hline
$3$ & $1$ & $2$ & $0$ & $0$ & $(5,-1)$ \\ \hline
$4$ & $3$ & $0$ & $0$ & $0$ & \scalebox{1}[.9]{$\left(\frac{71}{50},
 \frac{71}{50}\right)$} \\ \hline
$5$ & $1$ & $0$ & $2$ & $2$ & $(1,9)$ \\ \hline
$6$ & $1$ & $2$ & $2$ & $0$ & $(9,1)$ \\ \hline
\scalebox{1.2}[1]{\textcolor{dblue}{$\pmb{7}$}} & 
 \scalebox{1.2}[1]{\textcolor{dblue}{$\pmb{5}$}} & 
 \scalebox{1.2}[1]{\textcolor{dblue}{$\pmb{0}$}} & 
 \scalebox{1.2}[1]{\textcolor{dblue}{$\pmb{0}$}} & 
 \scalebox{1.2}[1]{\textcolor{dblue}{$\pmb{0}$}} &  
 \scalebox{1.2}[.8]{\textcolor{dblue}{$\pmb{\left(\frac{44}{31},
  \frac{44}{31}\right)}$}} \\ \hline
$8$ & $3$ & $0$ & $0$ & $2$ & \scalebox{1}[.9]{$\left(\frac{36}{25}, 
\frac{37}{25}\right)$} \\ \hline
$9$ & $3$ & $2$ & $0$ & $0$ & \scalebox{1}[.9]{$\left(\frac{37}{25},
 \frac{36}{25}\right)$} \\ \hline 
$10$ & $3$ & $0$ & $2$ & $2$ & \scalebox{1}[.9]{$\left(\frac{8}{5},
 \frac{5}{2}\right)$} \\ \hline 
$11$ & $3$ & $2$ & $2$ & $0$ & \scalebox{1}[.9]{$\left(\frac{5}{2},
 \frac{8}{5}\right)$} \\ \hline
$12$ & $3$ & $2$ & $0$ & $2$ & \scalebox{1}[.9]{$\left(\frac{7}{4},
 \frac{7}{4}\right)$} \\ \hline
$13$ & $3$ & $2$ & $2$ & $2$ & $(2,2)$ 
\end{tabular}}} 
\end{picture}

Moreover, we see that the number of relevant $\cA$-discriminant chambers 
is just $15$ and, as detailed in the table above, there are just $13$ 
possibilities for the quadruple consisting of the number of roots 
$H_{(a,b,3)}$ in each real quadrant. A representative set of 
$(a,b)$ also appears in the table above.  

\section*{Appendix: Additional Computational Details} 

\noindent
$\bullet$ The polynomial $A(a)$, mentioned in Theorem~\ref{thm:new} earlier, 
is the following univariate degree $36$ \linebreak polynomial with large 
integer coefficients: 

\noindent 
{\tiny 
\scalebox{.5}[.4]{\vbox{
\noindent 
$3210262679261654261076827566487887978952677993133439499132124128349163099562576434327795282069826801341279384389742180924193782
0493734613249461891575123312425537758926009708401449664848127075028378038127755754562092027528034556095987331999834544915663433$\\
$5191471292172711109872333557363189077486365006368398775222856555887438310160532137222984093986836707758049996466952978163430548
3859700433847406729339699486637419027417119167734317885395034944553549205323071048451865378745991325637107658891562004803963021$\\
$7916582410834573264574509022306809216807166720595650118477769710969358159991145000121900210603129391004020483285455074418130971
3789199968657912936915716304713300368271951125764741484911881789689710949909280469983530127882582940197639018652574316097165927$\\
$4851767125037143033108688125480023935918705794714022244000779076274084875949361686127070180554419254116094533471440257032398826
5887619665759127695813695117544671313029374859803482676720522828203332747369132717639329813267827737217645582653087186244634244$\\
$8315574538386240294871500440921685897086371334908759758652247944309462173537301927020295877377540959366871672811910466882363354
8823001738441866271417933408593742026621020258172531171037941687096683750391566890201196277854928605146109277932049948754347537$\\
$163973825423018619649916928$ \raisebox{-.4mm}{\scalebox{1.5}[1.5]{$\pmb{{a}^{36}}$}} }}}  

\vspace{.3mm} 
\noindent 
{\tiny
\scalebox{.5}[.4]{\vbox{
\noindent 
$-
324382459094727259401289824270892717683868251696388565175180972286962203032667917545195263656728566935541782620269589683359597661955961261887707201262224885857425871471737264699331797303883803968394698864138003263168127214087938952240674743611490633075$\\
$19850701989503904619332232029244387241205393201440694810935873640258616421750065659255245363626127872426609511971855627667197313393503943935070599664392437709509088502245387925686438289946280196755777266735568628027670302891895942455240380239348541007138$\\
$21876272802211117040616588914133481542941824682064399218207027410441020197053965792184501189998974137096525788481061263780544293186200272521540353534527194716541212710998494979768249596412873638118465722088755867101240903394314128286644743262349682863701$\\
$06810908678108796151682895256154423432438908740085658620066907509054080255602758882984730490251200403864036957003778027884691664507585182259273371004786586760987924835270705234099424031299758107216489492187764346506850031956602543757608036219422870974930$\\
$67732634931150836391329536335398087796917789145484149421322470928896660882181893396484055156757199839742907693001348935549347858623872619803040707831797500991801984143529665189546625034618206721153912607485436382175784883313631660224705891120685005054479$\\
$3814661669247778198153887911250960349921280$ \raisebox{-.4mm}{\scalebox{1.5}[1.5]{$\pmb{{a}^{35}}$}}}}} 

\vspace{.3mm} 
\noindent
{\tiny
\scalebox{.5}[.4]{\vbox{
\noindent
$+
 809152141556974868934399437864510021115111638199850452154855139204653323511857119330642547665300989945870794064603180964926598923609955752367605554217073862674439647229998022151077486688048143340165544515968827525195295696625130425168543995639165121020$\\
$98606475393111111131490372583540978794759332829389053424271462668445950957039961943297039024533686426650133906312683017778295240585842333997994436209588302327316697197347868613482078999514610307665779831182274397293156956076163981568149348319521175272240$\\
$69602848973638223671234481668708773638524066772208569233753719174414297140322972755096303120288902449542746379700091943328682290567508472122938525119335864004483482599813442688491358824874463785708689110046118567236108587553888056361225118909448435239976$\\
$85363525778275063882633985088861633045637568572418640655599771340966117316077177025548354465000420879366141444401769488743023506655227369811429070325743097357279584947389722894148243701640335386899918828074158154904040967644179194033943590484532588868290$\\
$00436208715484350281430958521701334038556520345798742867383297274251560260125128140024964357179172178332729541373333959360835526599891646278296395245685973824848173011934794033873673820617010584089578405333112326939610696542537137786715542314111808170564$\\
$17797935398406765335359304883741037662712347516293087232
$ \raisebox{-.4mm}{\scalebox{1.5}[1.5]{$\pmb{{a}^{34}}$}} }}} 

\vspace{.3mm} 
\noindent
{\tiny
\scalebox{.5}[.4]{\vbox{
\noindent                    
$-                  
7633904105907547477562892596982173156360832428474891638591142280223079449105428706756583122732808646651417259942685946940819433000247876966977809755154461172859923440135864704859309895115215001497405657006089337234621718687906130148468224016776868659001$\\
$62906136161520850585055229710327615895373633285822659093359999759457458710371936805179203414308428102888354148453197094948592874609129412360322573721241810773069917253603026737528739223841357795181341214528637763764305940244345764896329157108994824001430$\\
$70548398716960802523707154894569814066726931376442771382285497448127584490189190245615308138436464081129639494695367112967698419661378689208295770957612492642291249004300008218070862839972858752444485541075912627193961053233252525992817948665112645419206$\\
$97714451182288364625472250038956557359521193567225818687774227515392981494239906571136285437433240050519932117813368886344473003490946753446010255062985463701817020580378745389688272905603252364493421923789662825946820639090850813643810344723788384542422$\\
$37640404960379482370160647016010038981152338776026585752705491788644166852790244968712368553798179456718857334479594265732678181280271278287706151946939442189904418249561365473622345554420802851671377314618012276853164560682830889734794490277163604681751$\\
$88619460331232981055023512455744226796741203958642157994400284672
$ \raisebox{-.4mm}{\scalebox{1.5}[1.5]{$\pmb{{a}^{33}}$}}  }}} 

\vspace{.3mm}
\noindent
{\tiny
\scalebox{.5}[.4]{\vbox{
\noindent
$+
230093337646785477642337648916039584053600576061164203836479856553744690737907759762946847149790738111543069537433654480712083756545542015393822404537775938178647276684006103213601049812550632291587585500521757128044456791791517215010499352685682989940$\\
$76769779136079532083077298706673239288893430547919245882000690384552663378525724682080590367066390003253978188720604781829296579686441989050007989251954136554958219796638883651480190661818142503279238361895902370860815163087686521448383900360646605057609$\\
$78890826105870835196938719810241258004330760321024151375221737944614455609688046667623474015413579792295677660860947762730144238064287141468691697521027820387209035086917480673822337459248083041738559552350364778364668274510025762063736766654818143500548$\\
$00163917475787193035493871638239274568046577913648596997972933624911256666134380992044147745362362994072412007500079309208169709558009697853516590623277219165261650615418287068381138240915728417514547897616382409417472266723898018859911666862414860192900$\\
$65781223233196789579812681408173821793011757528848210543866309738400286586930910371108950292121712682531701995778219368135254780323667331349317282768889496490852436928966412532465587670544190393845777363265614996907735628907425526618557717660950995156409$\\
$03196750882212797881320616480430560035197835480219337487966281974268559360
$ \raisebox{-.4mm}{\scalebox{1.5}[1.5]{$\pmb{{a}^{32}}$}} }}}  

\vspace{.3mm}
\noindent
{\tiny
\scalebox{.5}[.4]{\vbox{
\noindent
$-
212729411186066712403244444117974895592083432277069807798374200460714586522389427190989578686255727457732429888010572012590521345474144744622584556492576226781530100021923781414167738836609830690169346293065701470602593975399056915288678831743090040902$\\
$70020651050460201278604825803439874527059691322207762978056438143083639595004608570770199659066173535049659236122214302019967547705913507396045630424189019124600912899235828450201587349041694868673303045690922153433225483182343034508287927359926758307184$\\
$29361359674851180729415678603257080492896008379132949177516954539414635873026733824163235429380301429922918108486141744501172967479336699601650583959048503224979357762373750131857798921399229764582782013822159119102048108993885078257776874724679329670527$\\
$81924424045644180358031655767360572539636934995589139380553928117046199828466329651137865327663353870159803103894468655021510005139643711558676606898177173470374664966362525756331378341194145254918829102850302193654213856412245627972572014978211011551079$\\
$42971492664535500221846353845203268223773030713567324081279321152468903757892108905525610936046750988562620372924200255772060372252704570517074102900494531105948403827430839082027728306072377340031970622545754084310727295005229059380041839699336147114384$\\
$437800392084871966566534271801493428455669602776363868762758543515101003400609792
$ \raisebox{-.4mm}{\scalebox{1.5}[1.5]{$\pmb{{a}^{31}}$}} }}} 

\vspace{.3mm}
\noindent
{\tiny
\scalebox{.5}[.4]{\vbox{
\noindent
$+
628018620374047481607618325611545680344707235799229696155423812514239761449075110170783506531390084215965838979354228071890586656911708186507816669827818512478078741060565628137293973341581177032784230330564433530949175141357121772941260444950583888878$\\
$90783312351874796109270117159553483104159972003025607101919561601072287273288673725283759701140671583184008249520956504939096395669176097937846118090492535370059617911825057134407410950540654195714308969703356980823547418136794241047558089797190713728233$\\
$07165423014948472826595737271408273273210576167362055961625803496626317922030045178543698931544654286249459700314168954419577466973483810581932153259237657528482306187287163571954343939161623202270833030795708044348408270024794978199520194884591519657767$\\
$13067147874875337728082898457354899679605086955129835877570869898172165988581141507554588934596146829355493859676128124153320363108582705024735401126254336905549276110459061079036134462431130137810947642418637029179064934509269840231797410082940494721441$\\
$92298397058066490989580056715382292024090400602132603024616738017899669827255998692185069205487103887183295314543365037290574173042387251986214358154229572962490004098554551362476305277204788493723857584057722437136648101388225469793239201807861161924870$\\
$585661067114926728027692869428726095287694765701441240820897441131522847186780954820608
$ \raisebox{-.4mm}{\scalebox{1.5}[1.5]{$\pmb{{a}^{30}}$}} }}} 

\vspace{.3mm}
\noindent
{\tiny
\scalebox{.5}[.4]{\vbox{
\noindent
$-
915407758866453063188890911462539032833543201092169934373155191912379892179001221867978554655467330356620179883427135691452186385604189301195701868746633249896052598734373799491337587983708683441044016892303581610425571027930082938471663532744578038920$\\
$88295629798957734398247219023045250400199519656224360337091817992474144417609997313586106700956905260624814704849084001226600490980417340859399521350656578812378777216593659669531832222399829228966675552215231462734888422371280496389101636246632093837311$\\
$27693591765274301969978210926081075620113807077118097243629041809971421119687493205936217041920030606365800095648959995770813198311981074582611083651457270808855721538530396485869671272041399503420704155226529489376898999994131181050204789115930377992075$\\
$07626247971815982277735761400528966410202135012205653721044427751445031297549841304213577908229552469542253792212766571545798492376295920584742610069032246439127771497279220371080776886793040345274447839566649781281378504616677139525009557678334109534182$\\$
17911723390135160940742771994168265314751533109319746363164188285304546358642246009423684676267699239386465362693060146901221728307674022964453287204940191691600441688690428461958200918520841322806908338675451451813118242731293718239885865276138103184961$\\$
531272189813624494629189928239685072708823008910033550430979881033637945078638201012106559488
$ \raisebox{-.4mm}{\scalebox{1.5}[1.5]{$\pmb{{a}^{29}}$}}  }}} 

\vspace{.3mm} 
\noindent
{\tiny
\scalebox{.5}[.4]{\vbox{
\noindent
$+
780138433306506414075713805749662707693006884233362859698964108552171024704319951216950407888580413540044770462874827736402572203338866493382909795770057313281695057666655139498005707153236746119876169080414157009753999848453450579299621881444993016261$\\$
99923166016283260298550698161448493489856903694026613594738070743052727900708893631746399525692223439768777025060264615361180274520191503273865156634415136029057772723460527625419320139970332158364877828270831737799754981660805818086027746874625075709812$\\$
69745314245067258139928207788609603851080000493409875236841688660791496588179700158648273846068928162043862291196946428642724420915724782964600895085264621176539667583040553369075754756340572866306080579549377527464644585721212205340021653966493960521182$\\$
94898134763287245851974937039023890452768685594171080878476969968951606937819452803187029242654564219966411726671498408026567683405447987527255907622146603216922510269909727939063800185218125103947840425969852273304030433719462909706120522987040242684385$\\$
75641102245219098461406820848229434065723184087354523273274540782560855454988704162655016359739973488384223041556024806666487818033349314444172882438366181723417476628201515132980675006427436357925941191803359919849980235085724426784073017064245867986812$\\$
373762589522520183742538081931837681061720832495156426218640890377219142924515010556543872498925568
$ \raisebox{-.4mm}{\scalebox{1.5}[1.5]{$\pmb{{a}^{28}}$}} }}} 

\vspace{.3mm}
\noindent
{\tiny
\scalebox{.5}[.4]{\vbox{
\noindent
$-
424438448163669191137581393673431036760776523543799874160396600899238736678206233281731106947147659961101743827036632007770640413021871089506012410760455082314323680746000465984688237866770212465014445181164100274311537490511116596589122860940739646239$\\$
84748505014750361359585260497629590946714813438726030857516476806357514220906454082139163146063492603537481285215737514617319567010671810077542793713164045761639690859095003045962550476920354356776058885890643110325903147617960604980885417720339647687003$\\$
79043133256949154390763697140154151602819477221599700940825472581463623072454163048708769209319349603855401587041018207946418616916079443263921949236647980016604770391313683581924350751895030608595548054097588807473496374025484851846439514156909323550885$\\$
32471203596588089684290854127360704163201683093327095966976928041763171286042291496186367145637527180508922325032542423356720442921555517776506992305131810465526355241448844446035458402585480809833645707732974860437197502430882898482820981266547686237666$\\$
43465463549357412695122320587732613511173170238475059846488320072652105700159389394782343845241096990990581017090309466022932922487821057562947263925128028273345125140315703620899122766709796487480596298205534428674965124119256491215855507543016530648158$\\$
503826822409525224219035729595205232523091304142022754194817730349898820358987390045673944784887925964800
$ \raisebox{-.4mm}{\scalebox{1.5}[1.5]{$\pmb{{a}^{27}}$}}  }}} 

\vspace{.3mm}
\noindent
{\tiny
\scalebox{.5}[.4]{\vbox{
\noindent
$+
153123902879208377756172776101133054643262097350809502848105165237269809861626063202514187635245948590018046623599772668874348828778621809923574347159362693158033413830987458982009585681996984951795902847862076269022478241675682157226018837118548066946$\\$
12265788801183554041601265568412687768079718601001308263313543612836730643658724087619651321362989318712759982292762191009132643370930965603870841818214749178676535896047927969912778581350157041067676728833461556456572205978056490604045862792243018035481$\\$
71742869177451392574484563505303887787077975029100108732148091383430201639519112550674102129311078988034010159915109982109736221365837443655964864105460358090755233611795203496985366077860058281087271632985054502628740262489085676302196457139284297458092$\\$
08326706662883860915886318055272568200578364045810886387703978926374756957072232679842586687367136688324954498601813210892496687551026869482859536112078334115018435221602226367590814494958463636266101620429115715947971176434197502040337319953364918493223$\\$
61744046705589399497256763150237057736402974189754780291083200698112196921853119511679374927529873942636225801764223582532061269838149658215936194979965071206949381719475243275518649788837728991407464577974205783891542303034214071805148798887398316727806$\\$
046971713748681300118979382735285320421531333031180243239918832289561103939011368108068136182419071152123019264
$ \raisebox{-.4mm}{\scalebox{1.5}[1.5]{$\pmb{{a}^{26}}$}} }}}  

\vspace{.3mm}
\noindent
{\tiny
\scalebox{.5}[.4]{\vbox{
\noindent
$-
357334252042429205326719380770071784873834467265740149494241579321165784692736685795715191918241499643346513385204838401847776297606883093556259375778335738186792870451195052639308201751961657831027766590833591260183367246633287505576578614629932824452$\\$
57151943055309523989669997892494282636671252370785694666358556267700935120994108893350674765352453230144284581592305957429943989880953482925813111995700138376451622514106506984445644244237679364023688875237589879136996274654818080784719748554772924268796$\\$
23608988466473259540269945279382946420813637137778645755248508411506035641432395940754285306655358018665126182177494275685701989742869419493158105575985400352599725055081238847301824724040200877170154343524205600267692439948390106232597979488770294545229$\\$
58461796607350270026307725804984223428685750596727785089795042971288998874056308744672537472267833935317970740435566391904901107209128399644592356362588876893762394690138343534615606218687326603614408987978660865672089943209614519736603191360403257124056$\\$
55874290251933149383933446078791636497848934481446954899500396716998869178788006251506250185576143661972120268860620528354398783338268240623238220484104945020822001785679919871803411954088438864427086770022113732376169606541040161684539744421377236745978$\\$
07315027489382912173684437382581332201859382318573587598790666471424793003993678361726221966480816458070613223800832
$ \raisebox{-.4mm}{\scalebox{1.5}[1.5]{$\pmb{{a}^{25}}$}}  }}} 

\vspace{.3mm}
\noindent
{\tiny
\scalebox{.5}[.4]{\vbox{
\noindent
$+
504539782645874758606053715808503571625726200484969610116767625750451834805320357066384259859238231447081376067561160270955552009146337145968314227037773888340252523544021546613391432314397962415411832857714350062204903902180623262062145170479195035452$\\$
23620757638478418196553883533702406174614714158523705806329325616352327284641924040468823818162002200694001468776351757647973547206106458353023095635269585155599351386055559922981103402992095049357984275902936359158934083826937063556474077338326344402602$\\$
56147650810677394088426351706984976446881749617669329156514167386774270039886723812157490445475310062244461347959609427277978319048531398822243979803490861356976655584307668482756991697353864292396496752198170950428279885597918126209572130854243147291395$\\$
84848632546721616885646194187235729428540988718939042918300530833551025672848594995260245022679736732426700719779886360151939125161706958692625909519870363783407225539800863379035954722767833690999891911281459023320101733099342552165962643140675370200966$\\$
95279366838837220569174101198257066456544192587834834251031332664496458229561649206859592538593221830774848921576945894699584415055804790726960106231492759410362596895389033790076443273262132877416737602946942154228789820180166778705183522028684594210826$\\$
2502171093132515153054772407053426083783064034889362978564370953528341645775098109730200908011806471766934912343398678528
$ \raisebox{-.4mm}{\scalebox{1.5}[1.5]{$\pmb{{a}^{24}}$}}  }}} 

\vspace{.3mm}
\noindent
{\tiny
\scalebox{.5}[.4]{\vbox{
\noindent
$-
447046047314625271617547405381716917060692343872903992338239579203715796423689272287549752753095036039400889397068060401831416071091640685577090336146802437757134945523228951892076694551569841395653802784684131706766099944713085595695869422809667471221$\\$
30257223803780883345715689577908830421545560019923117605526859056335900358916351418553211576309539804558673430536241681209522696287675032263535813211538785081269230045015296822521129806480909477218715904459639922252352707454348033823272837247060238345681$\\$
53667863506046270430391633718698959707619332991325387902289934637383550276616907379855804654171165626267609502444563103327394281576438079633821234602274208961446438072496119495734311234426873771675986039584728070873599158875646343800868952057422169285211$\\$
87033068284081767778334283596922962077864367526933717347446198949220575410116443022361746666952245577905148024779186090995933243694670812889092259683390865543720778802821072526738367813097749346866813427789100603017609121303949653185443374426629400978667$\\$
63476896771182799823809364753921484422700950156087406552123925275863018614000244146925119701768363337014921737493180104533192738096755192874101897442619651125336664372891410112061896003542172373246641894880561383213296251028496965224246686291279424648527$\\$
310913111549397119251926824430515928022842117766264364511134747838991086987479466958262829697926792111609188078745401044238336
$ \raisebox{-.4mm}{\scalebox{1.5}[1.5]{$\pmb{{a}^{23}}$}} }}} 

\vspace{.3mm}
\noindent
{\tiny
\scalebox{.5}[.4]{\vbox{
\noindent
$+
662092723221959325354008590078216748302916508192901678202998431888193521491316337470123022815222600129591542696871700834648991547796092585405486988793961844360403668323112225621506087443910087963345861259338512351407426125967471718765766781622774967979$\\$
71258001617963569099443190986109584786767425380482628591429766448789541784860330052598780564953242230949307463966960097539246666416252157192933651859506634850748626359026903428400007939004939729777425749956494865116428294184670280279016556973368040507865$\\$
69417219152180984590916653542929015778094022723565762077778929865117841566685402907946055440213254097988087064256472819463603222777926752335335991923412720094650484023534676305474221290759840790056230004712975945188514287896376153598503291343307939293513$\\$
40064984944868081120539692749330082703912053490425820937448734154652753814481524653283594876238316878455567852674240024233002964685160424468624206194615711756548595827598143278083700402178547518341409059986240116325708737986457177874026044294252404384192$\\$
26409238874101199524761350818576040406560547234691136897677395473328057146328473151881779369399655197775048766541167027562455410567410677381650126085554271742783900865279270779784811945281884021782159353047575375676622188916195895695968756124654742126656$\\$
59365630652948294313910681666378665516863770165265968757177967767767804781230891681561701210033658400584017897586928696933704794112 
\text{\raisebox{-.4mm}{\scalebox{1.5}[1.5]{$\pmb{{a}^{22}}$}}}$ }}}  

\vspace{.3mm}
\noindent
{\tiny
\scalebox{.5}[.4]{\vbox{
\noindent
$-
176980183801580805487541696896277288060543166557547036545397200476834759631854573270165708384080420067844500831259663050401012139245035717460638014324222975969206375483354167197509352661461808173504274585770087857908539608597845354813840302730175467475$\\$
40943838784277898677923010115946036079870712977474537622127012352756831171450373922552633817615852824218821081350335219601007370035589417173069954835470403320528698933809360968648170516386948827467414558802619470723330072645795974521626720895121099586966$\\$
13910443146431262245660509177056971182513121920990142431027683088901270405545286616580415767719654989105340585895016639771470851514848763245612956432855072267863465784583049941182759162996469351316564207778863803397732150655894137421834191227485973979168$\\$
75187192895598209858839022998752395497805821968016436313292791203616430736917307079732618160297179790308748519280091434671047452422156756602596744843772974729162317367075748482963291447293396763291205530032076503740960258564150457154640186841762617749034$\\$
84928842095182893651784112544629768697029697382668244573821858695902378703332052331734403647395579334135462158824896313457563035263925658285435827028479622089484155079504188276326360096793670053919850036211359651715463664731964734449181525988508364624804$\\$
79099383156220008420009952370542012933064766143343040044944587555967367096421572476298330088178639620210763127505252947572015608029511680
\text{\raisebox{-.4mm}{\scalebox{1.5}[1.5]{$\pmb{{a}^{21}}$}}}$ }}} 

\vspace{.3mm}
\noindent
{\tiny
\scalebox{.5}[.4]{\vbox{
\noindent
$+
195236540165784295066012139350991346500459039898959019714578567855984574970354886897301896569830420087361150579158111990524508666560811443030066119456398537360317876286938307802500353386289988683863754008876094862155486741419380529784685782486182979392$\\$
47394458863424467599178587603204397248844743961909805826588334294355231460596460596159547708480531953325672238809546133803682455648506284123823843351222307842539568990457749345303824848055486479173280619926084273105450078337868482168121076675991645113284$\\$
73203722754481034430896132259819860133616836970036903557642588244305286391006008042565495973372227513566314075152144457711698542332298802468945853133424938441966521921932192382419810197418482684065512878591105121748558477729782185649695520727835199878801$\\$
90477440661373394668168092370610484850010760055497216038456542229373098720714216822272779484129727927237121508466798700528301418682092927966297960100263736311791246273094868553352645428516205371012726897537007221895766917839993108149070812739309737139390$\\$
13154415854961098359576269767196566717219451764502139844592433604162511476928264847470465020575855372165551065360988478250939719959112768020213179501123449581123917404102208337385300812592355481478627844139582198194259934702149248429702822104578213954956$\\$
2284082138207027664764066182302802653441834040129258544280483372293867261260667218297035922270362520966129156614170129665950105967304999174144
\text{\raisebox{-.4mm}{\scalebox{1.5}[1.5]{$\pmb{{a}^{20}}$}}}$ }}} 

\vspace{.3mm}
\noindent
{\tiny
\scalebox{.5}[.4]{\vbox{
\noindent
$+
685845592426827827943055627803175073653972514358661620527947122765015299789494937899175285999146711043761049121463742123621006514608543743662680361569596110906236344887540345038247573079952319835605852167057161957744480042992662294433292711917110169019$\\$
04720560821995573723936318187586852855912597002586596933739657677663926667228677405028780107328351538894298271132200636304734988073474687507489688016498149126970274263903092031567559268160501568194153885002135697832369914702855413880655425971725129879604$\\$
25008960113482405145693118533820991461252576412618486747680217098575499911109471211628667445100166093904558968815037958272763847809772381551703669450623866295951064852323241963301900836958400818521454778958298665738004697136803604435964797170311180236401$\\$
47442154949549816907054503560220587546116378136958576999063605357864373897979434862404583660629412551877180867076422078069518662086597403726827638864395554133301053189863883555304457101081820337268055184855069341863781031291555240714346756589534102051063$\\$
37699683314445027146842672550117704009122121146328817813828147053076889811265135967838564616472580504239630152822716978720312978515125474788062876408113660835116030368319348473798934329281967840883939178304121101330199505483332726523329660390923589082876$\\$
29986714539872373938998661415447011222490344557578156375818949902742958695461537329112185957244859517552712325351174668237930728692009079104602112
\text{\raisebox{-.4mm}{\scalebox{1.5}[1.5]{$\pmb{{a}^{19}}$}}}$ }}} 

\vspace{.3mm}
\noindent
{\tiny
\scalebox{.5}[.4]{\vbox{
\noindent
$-
541547788081761030352948754930242960170711660619177161336162305026955406128804472053148547119224451154129268135459548617423052893010769642777109303442599872307067781813575892417155002104643972515389175472885243344060029604322690701056903137004709167985$\\$
24014937486879324156036621715425604602216825994635461160003521543524151497615958787739677976986951935674069377009585816600078866515127434546786642296351527370092046856047405898553220480457592104295426381520250594369442576704513809131188451730327945523506$\\$
64678380970732094749374455330293737600456947787704353955708817753226368328526438686154272304164612554816412375124671757288868854307922741107368359019664658558485432439487470248295523040225421563145976055131582909065410535269179008625327268719079533851925$\\$
61139665893516024429912719808128659447275429577222424099636159354388666516705123172025370615060932603822892749680840362775053167358155035864271418053215597036464987290952420254083525448330573285782275580847937019866539324686834127924789233142862667996111$\\$
12821393398642035366891832759983051723222691218860358058382438161773349992943315151645555959560834155658798993236771159204422097093346291183022246772445705768415593194030388670955908641639757918353664871011590061770071550138308597312193242816717240975384$\\$
4903725063658469914933498446331809483590675822469612528454399062628281049602185203387322744323687438794184017177638379262227794324396239598000770383872
\text{\raisebox{-.4mm}{\scalebox{1.5}[1.5]{$\pmb{{a}^{18}}$}}}$ }}} 

\vspace{.3mm}
\noindent
{\tiny
\scalebox{.5}[.4]{\vbox{
\noindent
$-
243814929807105462834684062983798282704408160958569118695348543479671192893599436152005528331450524098881703299035083914027584953888457782786589722797348719152513502515097540108329347942338104201232803105380592389084113163321937995483537984852691640239$\\$
10269507598154622488896472855525016374481110830187679057803343329385284538719389774874279147990275711361780510937537298889866935289918757703176982222223625245370555299712077329325590110481687969439796289560461093033861265370385947928243564767016658793907$\\$
09074943865297300888154412567763654568219184382252368072328974017615468861641537017638126069536390087037229750357160185417003880880367618251279526212093335458579735758902039583442467794705674069301551451584246885743965404408624794573379573741431401156356$\\$
39919062673255791457731349780904114246159824603423911180282062800347652521887891440074037068922830622304507304631770714470191402497139803814385385441338017251753185014394751373534688983067404267985996048745490110910197069781215014981859329176429074689675$\\$
39754830672216273641574784044913851974028838726519089356710149355274874405793615833349239540306804109027559908006769585431210984222929017038168862263400649409582371940122913407700686887097642480128488168996481384058716283352980118768669386602653644111964$\\$
2804040697357792768755453454210257377575462931428981674683331949746702698791999617585768961061875234644082786973320241072717452696140911868871130048559579136
\text{\raisebox{-.4mm}{\scalebox{1.5}[1.5]{$\pmb{{a}^{17}}$}}}$ }}} 

\vspace{.3mm}
\noindent
{\tiny
\scalebox{.5}[.4]{\vbox{
\noindent
$+
306817408884675823076007082516924231761562339451537116873127440540803324989939548345360250673518772281957678114670790847675588292494794092418423530307810307373082994702676656683189033439484064560598392158361952107791179852833540323463035513426870099606$\\$
67129547992926565683127793076269168447472189591143687883722544733216533851342874055530284830615640009548868868643427257115161206522238644797414752610750973852576332150503659452862746001383497401905055355879148148010777523102463867385737787642368498399706$\\$
73135159021927320198981185389560369999157767171031230943072492704503060695630853269035252530036312985986773858086488387238945441174319853327306215782690206604969927917962611707269188770684554977633054452319532990976699562471227270235557531429620055277042$\\$
94274551386056015987520282125941421037808476769211556248660486049486081417133028443452828134101394915321614561161425861800519706051290753617131502394823595052380290654453632822972224228674107582456133235929131593169779396693886446566655194712903404612058$\\$
94985595934107142130404875776346040274908419746923537537228831677769565407671499888801624896253357296695478338318916890017739888608157350927844561143814102780688726997894468273211392532607331420599828822569449285848795752968571199308671765648662296465584$\\$
49779042386126402138120667580148163670007313551900219583294150364003580210919111644911452404152884286982964786291851753154881403585525165997673031171218534825984
\text{\raisebox{-.4mm}{\scalebox{1.5}[1.5]{$\pmb{{a}^{16}}$}}}$ }}} 

\vspace{.3mm}
\noindent
{\tiny
\scalebox{.5}[.4]{\vbox{
\noindent
$-
201986459137987840042074460800467134095054583087261968275708636224411092028151531962817301127749024856080401544097641468335627194713930892403933476166937901917948096612678065143119389099505672176561128588763938113984786809390398804544755391334791709863$\\$
35910423495894894310920766158273367903882187825815675065204058713126925017373699302704163959309791889106517438352692673889281305946348716271466267050984559776372657571569945116081818801860469602959108161532400353144463835998010162045564543984154729185498$\\$
56446856104853652287641856783134349348446144104657798408243010680625101154786709965704247285083517435481226865851753017195055509831992360872101450146003997645115702457518080616895333546304122950624848059403678267435976550114893177193515603491461204227097$\\$
91518346421127109683743184773952920755884955262432559776522828077003469312488870084559192441608372510477722911081588959998406304574179244784226210468508822939056653353720462245655043494381411047608721783378438715433728759700234988862786614005717356616000$\\$
68095762646580921563661013060360172050988714030263391986007289937082494001676702511665108176525275793808128920661218480264917024261742511082271279254478882476342996944390238547349478414868882239661374430225233642719780771094503123629940726514230783754035$\\$
720199472293449494216701659490892244254583355136501025280050906887215646489958675970992917167917079217451837561373696710711099467426972313930467747593415951753674752
\text{\raisebox{-.4mm}{\scalebox{1.5}[1.5]{$\pmb{{a}^{15}}$}}}$ }}} 

\vspace{.3mm}
\noindent
{\tiny
\scalebox{.5}[.4]{\vbox{
\noindent
$+
172610880929873922906485600009516968145433300227146293052402913104542757262001951172278619207255059406450878809207448834876957866609823276722719952110237677126792208642079920667830396226130275773460404885835871033548539684379654713212450947204725155909$\\$
87876814813183906878413661848178280270084682266558696856868075584979577851137677404500108210743302721490706689138231551658910604060627762453858696524693024660531120595981718864648807187416171105643511194004893783644836862751987559808073633570320325470319$\\$
21458830592142745207342512085403784059372371542160296359903000255939660566344629313583267926025367430894462225120585130603675678625230578677539802675059572550243710299506810824024269155695782012106487537878163282794804648606048725755342737112558187560811$\\$
00732984342843409840420687450947894054784246082867561348790770041629193180201662956331342392725015551421908026708138822301830441139991653766555119846685612698151500952549448506299715572543328071439830278773475978476434443077743137626756314516342023908487$\\$
54087544485827650722771167158203454406435614297168313187523718648013322672700953829008124985053403526529573161926598971432614616946443003148291158236382803041413515042413976145477227729595541956088958045274732463366934444068394087405064669903816557322205$\\$
864092828891837651137061736256370521124845170546874118827762204236701044208372022438444385220420899260153437435744635940988351236479678339193704137823062657965146243072
\text{\raisebox{-.4mm}{\scalebox{1.5}[1.5]{$\pmb{{a}^{14}}$}}}$ }}} 

\vspace{.3mm}
\noindent
{\tiny
\scalebox{.5}[.4]{\vbox{
\noindent
$-
183353611419965929534156414102820331049041592029127680574843425318164390528270128360456484272247414469306721054212146603160937921750155318212613096014139613896119499110632717247525720313410639424342584022413138943723889650361575806250775884076601478366$\\$
56400133385654576966662029668624432551927733164687111276902509798083418150454515763063403399078617067883232919472616624189901406863382187577616715257412970516741698625171198211859419420809959737951738261910303242118957443557586605355209696328764508281862$\\$
74203460878451828772679638470204418328179246015859470991418727600726012329156125403023162413921588251231251864308003558843159925655801399884737470283767759657193850175507244637763341848392464725273104905700494108045907408088634814157947971969864079825689$\\$
70644399233298779945275603330489399493143351267935439971500756559447461073674668053674118069977320100631807787593463520332382686279420641369583369223091482368099417843909182800236362958235015779776069999112195717297940149920891646616981308239057223055581$\\$
19734564478509099709335325998869849057514978362285684337511741709837391369873152050535252078361696972111457534522711169139176463915502852668286292438940429602111063835910604001586277581208451959835328669651990136402086868962292096884520957778702809626596$\\$
870640736907700839782301365917916735256961745899875799888452523272368004381133589132835229204242180495075650644728947440226923634946083829345065679906440000000000000000000
\text{\raisebox{-.4mm}{\scalebox{1.5}[1.5]{$\pmb{{a}^{13}}$}}}$ }}} 

\vspace{.3mm}
\noindent
{\tiny
\scalebox{.5}[.4]{\vbox{
\noindent
$-
402862236608300847658017381141267392180500678905652834487405657002357780084224628754867167494461700730844943343447837931191152103444482926001397135851642694613006381466460265289330920160877584652261389728284230373473157197797837655855626670762015342117$\\$
00419147041478603852662480920021498972194312296806711961423454721367611272492406326362734829692664387371014679322868201379378705686578504840420491504421302871465487065821693655585284661802521163249015820514854068411840121165058505587729797798254887311085$\\$
30930366046060131029730522084755041525419385938944240261835963899187777955226460200057349362230057428014491969897801186064021458410249919829374385063875622134923734427291154580493963201662590089918252322152547487325213252324829857104705407925693504401153$\\$
24661073680850080251145060408202413198352971999410603090890566952291838157161169679370385243490159254600044660108116764778001206673003881214653483719341394171575019801512947066128952965283945302568373485203859426639341000019993850156087292897750194903662$\\$
20976710625119502580641803397907747424844258041046132455162285477389130706522939126914245702054221978861282976143768008337988327141357763566756418617988587628177282608254854377724747764070113128754518091216036096944062251607143904153031825111041001154714$\\$
57530338782292214418773591978176061281086402264734178420095841875372925791694902003254899455603356280826728540924935313074974186678049934059713745117187500000000000000000000
\text{\raisebox{-.4mm}{\scalebox{1.5}[1.5]{$\pmb{{a}^{12}}$}}}$ }}} 

\vspace{.3mm}
\noindent
{\tiny
\scalebox{.5}[.4]{\vbox{
\noindent
$-
133288041124691093608755297882318205780381119187982572593919290599427767881937245240768149088284359365561970581906131634939972234384931899076851214018728646763297190489302561082404819304157116322715072218027259958452802943853681038464284262279140375272$\\$
63570980121739493976617489761199341092282168593826246820577960244044228571608361998843091323673458585667000211518887637208072877413223793184691611195563427046434643464299952450172177572502206160877056976058938610115705297112616472240477138365815306048691$\\$
75712705430718039960605756006554790075708913393554163606371586646567697995669215493151531233760107402222879370195125137437050302317045488917980793887580655317952563436830843275597186551874880421663245736806720479081037756883874915472879775082519457767643$\\$
32311821658383343920283178331745861657323702605710409110295871135405033500373700821708308225163348056306568021904330674209710539794607815137874203407429228944405600882946356748384461697111798001559109864008312150228803205136545002000983115300959316686488$\\$
22364293900407960202746003140473042485607841327191258516461379822136676114640503610120418497775911278875893292275225815372594123540925050929760203604261939094184084575118824098815306576867436663243220349593375498296121384634368731571433295854070495225429$\\$
8283294861865519184367604423351585722426293778640732453416386327052558062873300012370928768157531000215104567902613750997597996672539011342450976371765136718750000000000000000
\text{\raisebox{-.4mm}{\scalebox{1.5}[1.5]{$\pmb{{a}^{11}}$}}}$ }}} 

\vspace{.3mm}
\noindent
{\tiny
\scalebox{.5}[.4]{\vbox{
\noindent
$-
175898906443417591139496235172568662170289523794505987485565582929341686764436373198340886957514988774587731765772676970626402128542972337935454688682744898319925828190591661160555317776499801939896516485382933360140227201204967364321110104445391806572$\\$
62453275488127933543731910340802019804873438733069281978931859906113199602651441619980193205081552505334347058229918104642608681597596668723203428641780572811197518086155787972309988895774379543240845660912935625709038527279570147751188440281307221342190$\\$
02485062742668498101005228267341254628522832201054769311197575152552005206823946545316898790140091645544290151657716759169280918239312538108274210373542928170898061853474212032197503088519851263137015568617321922252906789330757614887975450370606937248472$\\$
46862749770153664120982845364945554319544554454897588671230510319311385633451873037645398708266640127015966140852398555977500636579887199840780722736565851851414659752497793625527851406809800166991615348458901949151830145724610306584896448109339012654317$\\$
79577160842984879875827940626933594344054650228950411019689433236604984948204237852572309062240571542089002992128614108133765952120313257737807369563111329456941235532659291011330287522956188410868778637628318715654761448744861440218280626884402042846642$\\$
895373427573397946769795063175426793801418775027365476689572527133898916375768379585808654921863121714572482185864573441674840736936857865657657384872436523437500000000000000000
\text{\raisebox{-.4mm}{\scalebox{1.5}[1.5]{$\pmb{{a}^{10}}$}}}$ }}} 

\vspace{.3mm}
\noindent
{\tiny
\scalebox{.5}[.4]{\vbox{
\noindent
$-
751460619621329922194735459499945276224970270884574325031018542876532965232957953538863244937856181495534826992145695974154271116570000878374326911336565148354453413971476801715653960534669780621411939751671815003608178527404671704624321419556759120154$\\$
06463009798565292413161746683421971542671072645781387826547931725891921062008789985228155432300097079156463797863443476951097196361793964310934211112342558253963933166886261552315018108296634351550263418966815706172688201866911148609117430945720131761612$\\$
53332527160249693407395232412074263819711023881466461627528191995295486860899410963057708870120019974023736428682726673784465081738995380021490163585212581368722013365241384349701315494420501175457925070473673234686468847588860590963128970126077364698701$\\$
49913414868793504563276039721983787925242210920529679566267702189774653358107569170665370989594966869856258656333313493947965453447063611304580287347672844043693859299679211588125739073721659444468035622240194214927051715565712270687456081518089960258347$\\$
46612171417123961542183384924980786520007955969217826519989392156610049802386449565761023254344753019860290749033428226151901630300634172400330365236948331109596540799766687684937730404073918072318161592722293847486176109641315676644784956687902934453409$\\$
3289346927618962035432842534058940890719212156502061364177638343158186863067797003099353114782488432446111823769770061444677605733133418652869295328855514526367187500000000000000
\text{\raisebox{-.4mm}{\scalebox{1.5}[1.5]{$\pmb{{a}^{9}}$}}}$ }}} 

\vspace{.3mm}
\noindent
{\tiny
\scalebox{.5}[.4]{\vbox{
\noindent
$-
301223266949496012259528953721912737149523065983550957958464665594194612079603387066656897886729754170905492453326458644789927556745397103464115197072803793894001785497812483506668073567035470163010600209189783663530329775750137246850958073381527834878$\\$
04798023717098734531361052815664755875174051530324417831134109586310502827290662717459179115466181420973973876397947863014615191230692158080653906716090231370432931670650473167717939864828571276419604059016835433882534901167458361023769863189117400320220$\\$
66920227854834291326105441126619092794630550660300338609852983385680985697107899535168537108645284157035420126443219377401469553945731639132888163674652893372897566870290853969305990647293356483586084427306189548376514820866117813858448619988319499998777$\\$
09722449910009160816787509620296910664465037318745570806228593922646145103871675102008098109157936083113683204394271944820457741118268270009687749094703976687098269558270001830234215970250397356408270548613298684595694379944625775640436220232025924686875$\\$
77461091654121935713790918386255954065514672825761837542656916453895178328616461171307539349560917845041614166122779686498797286197634717839264723860127096759424334725701494267687606383038485090223328955528811327940360653861911149742023311977758604251186$\\$
231884405854015718528305869020579483402812255715570108502574544858285027442371129510324804851391056279818089741709755040262102773751529127821413567289710044860839843750000000000000
\text{\raisebox{-.4mm}{\scalebox{1.5}[1.5]{$\pmb{{a}^{8}}$}}}$ }}} 

\vspace{.3mm}
\noindent
{\tiny
\scalebox{.5}[.4]{\vbox{
\noindent
$-
161673003569660507516078608617499261359547641896799721858567522018366961215699423009766416945784760546085730389929160084558638924191582135732321958820211318222874565581806187848951418006714183267684698521492844901621244151177594936392164478316851207644$\\$
56309600181583319726303733781194912533406804464984002401822612518059156533912437713763397378504137544185365786566848565203774295717742487780234632799028245623850341687021120169999429637615468336814254617662685787172223852192905545263280960851232852221698$\\$
76863344033098945116402405513306113904617046819223900371117319857464278069982584535909638315216508064846669583599247899783494442421815628485510930904597483569306852921202360910288769104485849702250861222705113043494526145731783404463469289150203248720713$\\$
40455506634996326240855826664085868480733169485050464427507789196830766790910924618856361849449010196036497960996347717623694888325383592812093814306380395502417132564487034583366228130249372228725301760115056807958719663032640946348504983821522067743227$\\$
51279647498370638455995069935688472090071076830719027393582229814827319897734987465193914270824927838065105256026163061356508476344390544985110713659254098572584766413407571667619727641659325741114552255817869055319642923340681756957768548788264508453330$\\$
94022707776857225535683319104577391657791240150272095044405773641206881066385209286322774335901691474141897452798356175524954835671324199219966999407915864139795303344726562500000000
\text{\raisebox{-.4mm}{\scalebox{1.5}[1.5]{$\pmb{{a}^{7}}$}}}$ }}} 

\vspace{.3mm}
\noindent
{\tiny
\scalebox{.5}[.4]{\vbox{
\noindent
$-
819040771786514382512549021127206599431608106658755605746147147533720135858415545424766474498997507112014687306152774460183585754189009613311235577639801250094387362107082086121903106603585782912252290748871219807033217248626167904142748977874498145632$\\$
75301902472346547327204796306019934919189656664704906826812996703859196906745681570651998824247666582271553184285067771202708274384102755071113761211149334839071609001714006531263816220633664194585873777229978142642291607727136766895050990000403307675406$\\$
35281046743577671170767977383845197506059177148282340650576738246044511526377607789422877236709606626638356297435404264264551210405791581927447187538018551398771531528484988150460797728660806976710005142088692931993916207283406100821411318047784226247455$\\$
57895538746464945074586679907702681183529552606863093845697483351105607715003532194051778382039111812561109795799237263975995477537622293498264903015110450446616767820035369793498383217307297815437110123756729221931691794611432102591360034335172021369868$\\$
67120119560434926399123403114257935298018688932551800165815283681089107035420163727826423255122431872762679625602594202454624105304681064822858634255252086799948200158398300047814575002440585529511445275129481793956897144769829907648324031429748149144400$\\$
88011893349936616652164614656949983683003999724225811869129237678075862279587584211740996974808319558336975225260774649078227955493070738812733111444686073809862136840820312500000000
\text{\raisebox{-.4mm}{\scalebox{1.5}[1.5]{$\pmb{{a}^{6}}$}}}$ }}} 

\vspace{.3mm}
\noindent
{\tiny
\scalebox{.5}[.4]{\vbox{
\noindent
$-
909912155901077108994105894597892630304681316729032252958158112535143736524409870327909480101867394240506873729189755727000916054736564801221088616262563092728849780146973183341819529852326957486198273576374710211182048609997898589278581671519095824673$\\$
53831582958566593672955597028043908263509118360368700373726323952533245720388534410446046259766128127305002905869761270856178324702636804423117763839270846704035778728787963500524997493395993432666640430503771694880256192312862796110454108865385173163794$\\$
15299577754114590142623325785057101381257459783116263992230341276775487001167166682975043805664998132875952934820272435203853875083929730169569271442479041092080626406116047893493207779709014209410988742478547470638418716658004251930178027061516675790308$\\$
36038097334228813505576841404763612570188302292641196238522021789499777466122344112048032443303769295411559034498774472557148301771105478215845681609472363920295743019023661334969336284387968924825706000213836973627409042543906257148998061478163993542830$\\$
80604376634594062594119811337506726054718579667507636107403449350427815698192579534488936983574879812966823960952162688837325356866448720471774142341378927683361886047882997602633961012198382986237122027387063045111389061982968190128413256411776769944805$\\$
8794161354478663449291226767797038998032442135074046931778364913914296215038342722536098178059898784346193088550942074778964055822029737951445182986276449810247868299484252929687500000
\text{\raisebox{-.4mm}{\scalebox{1.5}[1.5]{$\pmb{{a}^{5}}$}}}$ }}} 

\vspace{.3mm}
\noindent
{\tiny
\scalebox{.5}[.4]{\vbox{
\noindent
$-
228913584088828724680063773747742202708086215418782900895234547137267845288606130600130937185975217275131946459170167205946679955724423567391187645452273948025383356188884630788521282813963349729188864392712824416350736507133495344090805969320985955729$\\$
93615085028894258379297368261683765978802739244082865835498209053832836411577471566993649100885194484773544428054426928622813748709787437824928547170573922931924347381333140509624344004251861406031981287506005610527242757383478420242845309587695934129526$\\$
15669086155344223583873286970976319677613113457231082961107795830228078987520176447671479591925384190143687944494806973465562331964311632884535126874781069424354266652496012660269155703062839528838971615274620923193473607005205671196255862192310263968773$\\$
92102876945660264876742210856829836316317719255926585938523054454922632214601062701725171262972764957351579880369115793942140043055989140248418245332516634480864947600744895747024871813138223979847922618967718247665535467863311065831799205149265245394435$\\$
59879166289519479236672753442811748546485447701650991885066156075572530598482553695273615277903784034791373986065535194365484371708443724068917972172114980957275280338692387743211267990279766140382683266379705691852668659485402194687639225063260761028608$\\$
187213953769939607389141438846725760628097878489352637087633507680615508902117328996631359118616120131580459186048077880484261724468637600468992610025509293336654081940650939941406250000
\text{\raisebox{-.4mm}{\scalebox{1.5}[1.5]{$\pmb{{a}^{4}}$}}}$ }}} 

\vspace{.3mm}
\noindent
{\tiny
\scalebox{.5}[.4]{\vbox{
\noindent
$+
170275607517806808374336639679337638278592814766761658208945155964114586892790524505281621271276718694306383163901731698457549758872696130264190941557109414638183246205162068322565323918158432951913810465054038765848656263962531842581221924035659521446$\\$
95794908697464915795330694172958932113795374992993268703968704712622448739604762879282623963327656851888186752192998856019644734232039262944356061149735364297202673547764650201601901449793377010967128699759451066735556638384017254389081971604321816825189$\\$
88489053037099101916818055402135648590887195084895203819179113744452285884191134931206926568693480040857056563362973903688322038970374331680758621418151388414520917260839359783219427292348372788007418582014725799774535934661800605306556080757529300312485$\\$
46091494746188030608100536195295464506068501296682255451262229749433140165164370494859792429281375465487683007557103167734817582819918219467719752766335127460721913712237479185989462494716144155945931974865422081231714545887236090053134401776782112306369$\\$
40270092561833037390413696440742253529631146311782568861505703000774350094570766086469269426990602816971052702056688252498085350353847811511855926176985712385067947699366726308372505895816612516855468706189047885713130688803804479364020123997854028883288$\\$
73929526065143121897237654379211654698962509824000157509881941198002658276580071169893727950472794414381659337853900273173308740596348377444862608361131339762550851446576416492462158203125
\text{\raisebox{-.4mm}{\scalebox{1.5}[1.5]{$\pmb{{a}^{3}}$}}}$ }}} 

\vspace{.3mm}
\noindent
{\tiny
\scalebox{.5}[.4]{\vbox{
\noindent
$-
303576296659867869289546302797616493751069572946291511163377571493368206540527991141377097229882570216635876838404567961949205927244048158332225646826696187972816792961290276258130792780803662816773799048972291239849940000821417062901039288962360186683$\\$
24017023581883728680467154674726604706437996580597140313468619387263413222701202568551807176530555658593725081845245765369466542156302839247218939650900872497766729105762610386629421163672619247498839663758848141087967386940393233086574415635939214093371$\\$
99700584755189368893350769009979709339789300286518418606097415660316844715642885056089255572902438471233919561039090676570864508109152091959059413698559557776779447620911666165928813348568747245389847194177450363852118346819223872385043147664876150675206$\\$
26155064153458403096105148341934235840394805160940894882543031881804436995895766402344060886469411335124939557959193867034118785960095230180471937348241681991485134091715861731416439138159908961218288073905877713124405959503534438367986290063085245063186$\\$
94910647103492889536553156209948461205487000767243414846840708426124858247103048693757181559598137461757924399388990875721709021497228411795077775113976723734684546389073788021276177000849949317132747230045175673033279112262806802282209057428626607279248$\\$
431820563227772108566824730715198187436736049939041989577811135253783450089869614831562713765438536243442312740994128213749826154509964537813338641969451714430761057883501052856445312500000
\text{\raisebox{-.4mm}{\scalebox{1.5}[1.5]{$\pmb{{a}^{2}}$}}}$ }}} 

\vspace{.3mm}
\noindent
{\tiny
\scalebox{.5}[.4]{\vbox{
\noindent
$+
210947679211023092927923111091195356259466202354606220621039183162714832051334900548448288198879545474550114682564875952265427063251777729817538512148056211720902695412506717130119041880612120542871383957058165156511748214818927592842534353311062437718$\\$
66971840442705665687835554909576275914454007608763077516864676452027199920032473073108771134299511786852449529421568766550491965330810472526650734311969205789562074610077286829472502175829811343067777214803339543985204853134259942854316078787752381072416$\\$
08459341104293057691194040469065741374773384741412789421654767380600544630832079193371142173401308486159769206994865828471044608759972574220050490103151760430656565786404430911706258603799716458850755391853128940063190312336902961633448194476498376060955$\\$
80749326703467866868229039760857133469830934588291367582611921639384347285579993308803774398652389182986653095266498943109644926847837704326396822630818464660199651983261685954464183405391467737144038155710695336495385414770255987697960557859835444149458$\\$
24153096694553007171317791389330458970099217189980986517131239620909986509179274958968340469750601132200753236314758109802590901511608326882512899496163188477553292373495065211933424982540659311376618339984673256876435278793760010789303766517386316357261$\\$
8582747793373142087465718462445000331361301749890778845617989361133200072166927663727156145425685523041006115294849130014447832956858985842003156863100343798578251153230667114257812500000000
\text{\raisebox{-.4mm}{\scalebox{1.5}[1.5]{$\pmb{a}$}}}$ }}} 

\vspace{.3mm}
\noindent
{\tiny
\scalebox{.5}[.4]{\vbox{
\noindent
$-
507843069055560760365941682812940423454044079836693079407832722631986588931262081035917411613035783088918143467414613463557170220427832473299897524119199681500295699390448987791656198187526652520311916835270111066971350664339163915579346358290752382022$\\$
91556123443187746844345153885047771473075245021073579674606501264330176694073073942553656913391122963257996044281078530660916992275407884012232550951336962175015325787457036813258893743558563644723256590415137859948899025019355611394340921928106974650343$\\$
64113684984328115274456982405628236376785073983543026064479007951856610912992524739377470836095406234156608515337374918866861294384330788684280463230841923361044394236585303128105665367797331673960863638441433612846292096665279898573235105780607973387532$\\$
71956866900690608872165346245784713273899319170572197839472637199228861276100271077183590478158120140151203335583006555051672450609676593523663009844010897713529255746985249333757260465844976813649864480697998028344752252014897838377083388425601632788719$\\$
89125483411581637015567792025962574386427527624912984094548712276772097041724991176608500017136716961037478788881929727838386338892162104869116300107184649515906731194773738462410873969061035486721486412588065628663508063028366883426485108889941041090486$\\$
9599282867239514127558104263889813304619189995806503190188771758815556004388692589767155179187421837675746339145623495112807842754264198159575016688904725015163421630859375000000000000000000$ }}}    

\noindent 
We point out that there are $254$ digits per line, the largest coefficient 
of $A$ has $1460$ digits, and $A$ is irreducible in $\Z[a]$ (verified via 
{\tt Maple}). 

\medskip 
\noindent
$\bullet$ A bit later, in the proof of Theorem~\ref{thm:new}, we used the 
critical points method to study the complement of a reduced discriminant 
variety in $(\Rs)^2$. 
The locations (to $10$ decimal places), and dispositions, of 
the underlying vertical lines $L_i$ and critical values are tabulated below: 

\medskip 
\noindent
\mbox{}
\hfill
\begin{tabular}{c|l} 
$-$1.8562718399 & isolated point (on line $x=y$), root of $C^*$\\ 
\hline
$-$1.1581041767 & \scalebox{.85}[1]{isolated point ($y$ coordinate equal to $x$ coordinate of 
next isolated point), root of $A^*$}\\ 
\hline
0             & the $b$-axis\\
\hline
1.24487176148 & $a$-axis intersection \\ 
\hline
1.41544129863 & cusp (appearing in Figure 4.7), root of $B^*$  \\
\hline
1.41666026637 & cusp (appearing in Figure 4.7), root of $B^*$ \\
\hline
1.41767594900 & node (appearing in Figure 4.7), root of $C^*$\\
\hline
1.41790510558 & node (appearing in Figure 4.7), root of $A^*$ \\
\hline
1.41821476967 & node (appearing in Figure 4.7), root of $A^*$ \\
\hline
1.41951679775 & cusp (appearing in Figure 4.7), root of $D^*$ \\
\hline
1.43683087662 & node, root of $A^*$\\
\hline
1.47813022442 & node, root of $A^*$ \\
\hline
1.48488178680 & node, root of $C^*$ \\
\hline
1.59316011321 & node, root of $A^*$ \\
\hline
1.60149022139 & node, root of $A^*$  \\
\hline
1.92733319557 & \scalebox{.85}[1]{isolated point, ($y$ coordinate 
equal to $x$ coordinate of preceding isolated point), root of $A^*$}\\
\hline
2.45494131563 & node, root of $A^*$ \\
\hline
2.47089273858 & node, root of $A^*$ 
\end{tabular} \hfill\mbox{}

\medskip 
\noindent
$\bullet$ 
Finding the the aforementioned critical values took 
just over 30 minutes, using univariate \linebreak resultants within 
{\tt Maple 10}, 
on Rojas' dual 2.2Ghz Opteron Linux workstation with 4Gb \linebreak memory. 
In particular, the locations 
were certified via use of the {\tt realroot} command of {\tt Maple}, which 
uses Sturm-Habicht Sequences \cite{sturm,habicht,marie,lickroy} to find 
intervals (with 
arbitrarily small size and rational endpoints) containing exactly the 
real algebraic numbers needed. 

\medskip 
\noindent
$\bullet$ 
Finding the location of the critical {\bf points}, and thus 
certifying the isolated points, took close to $5$ days. 
For the latter computation, we used the {\tt gsolve} command, 
and took advantage of the fact that the underlying reduced 
discriminant is a composition of a low degree polynomial with 
low degree monomials (cf.\ Remark \ref{rem:comp}).  

\medskip 
\noindent 
$\bullet$ 
The aforementioned refinement, $T^*$, of the reduced $\cA$-discriminant 
complement turns out to have exactly $125$ 
connected components. This was derived by picking rational 
numbers $a_i$ interlacing the locations of the vertical lines $L_i$, and then 
computing isolating intervals for the real roots of $R(a_i,b)$.  
Finding representative points within each such component then 
simply amounted to picking rational numbers $b^{(i)}_j$ interlacing 
these roots. In other words, the representative points are of the 
form $(a_i,b^{(i)}_j)$, and each such point was fed into 
a {\tt gsolve} computation to find a {\bf rational univariate reduction 
(RUR)} for $H_{\left(a_i,b^{(i)}_j,3\right)}$. Each such RUR was then fed into 
{\tt Maple's realroot} command, and in this way we found the number of 
real roots of every $H_{(a,b,3)}$ in the reduced discriminant complement. 
One can definitely see in hindsight --- after the proof of Theorem 
\ref{thm:prob} in Section \ref{sec:prob} --- that there was 
redundancy in our particular family of chosen $H_{(a,b,3)}$.  

\medskip 
\noindent
$\bullet$ 
The approximations to the $5$ isolated roots of $H_{(44/31,44/31,3)}$ 
in the positive quadrant were computed by Jan Verschelde (via 
his software package {\tt PHC-pack} \cite{jan}) during the October 2005
Midwest Algebraic Geometry Conference at Notre Dame University. 
The computation took a fraction of a second on a standard lap-top computer.  

\medskip 

\noindent
$\bullet$ For any square matrix 
$M$, let $\sigma(M)$ denote its largest singular value, i.e., 
the positive square root of the largest eigenvalue of $M^TM$.
The computation of the $\alpha$-invariant for the Haas family can 
be greatly simplified by the following observation, which 
follows from a routine calculation by hand:  

\vspace{-.3cm} 
\begin{prop}
\label{prop:alpha} 
\scalebox{.88}[1]{
Let $z\!=\!(z_1,z_2)\!\in\!\C^2$, 
$M_{-1}(z)\!:=\!\begin{bmatrix} 6z^5_1 & \frac{132}{31}z^2_2-1\\
\frac{132}{31}z^2_1-1 & 6z^5_2 \end{bmatrix}^{-1}$}, and 
for any $k\!\in\!\{2,\ldots,6\}$ 
\linebreak
set\footnote{When 
$i\not\in\{0,\ldots,j\}$, we define 
$\begin{pmatrix}j\\ i\end{pmatrix}=0$.} 
\scalebox{.88}{
$M_k(z)\!:=\!\text{\scalebox{.7}[.7]{$\begin{bmatrix}
\begin{pmatrix}6\\ k\end{pmatrix}z^{6-k}_1 & \frac{44}{31}\begin{pmatrix}3\\ k
\end{pmatrix}z^{3-k}_2\\
\frac{44}{31}\begin{pmatrix}3\\ k
\end{pmatrix}z^{3-k}_1 & \begin{pmatrix}6\\ k\end{pmatrix}z^{6-k}_2 
\end{bmatrix}$}}$.}
Then 
$\alpha\!\left(H_{\left(\frac{44}{31},\frac{44}{31},3\right)},z\right)\!<\!.03$ 
is implied by the truth of\\ 
\mbox{}\hfill$\sigma\left(M_{-1}(z)M_k(z)\right)\!<\!\left(\frac{.03}{\left|
M_{-1}(z)H_{(44/31,44/31,3)}(z)\right|}\right)^{k-1}$ for all 
$k\!\in\!\{2,\ldots,6\}$. \qed\hfill\mbox{} 
\end{prop} 

\vspace{-.2cm} 
\noindent 
Note that the right-hand norm is a vector norm. 
The above simplification allows one to apply 
Theorem~\ref{thm:alpha} using just rational operations, 
after some minor observations involving characteristic polynomials of 
$2\times 2$ matrices. Using 
{\tt Maple's} arbitrary precision arithmetic, one can then 
easily check that just $6$ digits of accuracy (for the $5$ 
putative approximate roots preceding Remark 
\ref{rem:grobner}) suffices 
to yield $\alpha$ values below the critical threshold of $0.03$.   

\medskip 
\noindent
$\bullet$ A skeptical reader may certainly doubt the correctness 
of the underlying implementations of Gr\"obner bases, Sturm-Habicht 
sequences, and arbitrary precision arithmetic within {\tt Maple} --- 
but hopefully not all at once. This was one motivation behind our use of 
Alpha Theory. 

\medskip
\noindent 
$\bullet$ The locations of the 
vertical lines $L_1,\ldots,L_{18}$ --- which were central in our proof of 
Theorems~\ref{thm:new} and \ref{thm:prob} ---  were independently verified by Bernard Mourrain (via 
the INRIA Sophia-Antipolis \linebreak software package {\tt synaps}) and 
Fabrice Roullier 
(via the INRIA Rocquencourt software package {\tt salsa}), 
shortly after the September 2006 IMA Tutorial on
Algebraic Geometric Methods in Engineering. 

\section{Acknowledgements} 
We would like to thank Thierry Zell for his quick and detailed 
answers to questions on Pfaffian functions; Frederic Bihan and 
Frank Sottile for informing the authors of \cite{bs}; Bernard 
Mourrain, Fabrice Roullier, and Jan Verschelde for additional independent 
verification of some of our underlying calculations; Michel Coste for 
pointing out important earlier work on the finiteness of topological types of 
fewnomial functions and zero sets \cite{coste}; and the 
anonymous referee for valuable suggestions that clarified and improved 
this paper. We would also like 
to thank IMA (the Institute for Mathematics and its Applications, 
at the University of Minnesotta, Minneapolis), for their hospitality 
during the completion of this work. 
The second author is also grateful for the support of M.\ Danny Rintoul III 
during the early stages of this work, while he was visiting the Computational 
Biology Group at Sandia National Laboratories. 

Finally, we would like to wish Professor Khovanskii a very happy 60$^\thth$ 
birthday. We have long admired his work and his tremendous talent for 
turning even the most difficult ideas into simple and beautiful explanations. 

\medskip
\noindent
{\bf Note Added in Proof:} {\em In July 2006, Andrew Niles, an NSF sponsored 
undergraduate student (DMS-0552610) 
in the second author's 2006 REU class at Texas A\&M 
University, found a $2\times 2$ real polynomial system --- consisting of a 
degree $6$ trinomial and a degree $141$ tetranomial --- with exactly $7$ 
isolated roots in $\R^2_+$. 

Also, around the same time Mart\'\i{}n Avenda\~no, 
a Ph.D. student  at the University of Buenos 
Aires, found a new upper bound for the number of real intersections of
a line and an $m$-nomial \cite{krick}. By  \cite[Sec.\ 3, Lemma 1]{tri},
Avenda\~no's result immediately implies  
that $2\times 2$ real fewnomial systems of type $(3,m)$ never have 
more than $6m-1$ isolated roots in $\R^2_+$.  
The latter bound is a serious improvement over the best previous bound of 
$2^m-2$ \cite[Thm.\ 1 (a)]{tri}. \dia }  

\bibliographystyle{acm}

\begin{thebibliography}{A}

\bibitem[BV06]{basu} Basu, Saugata and Vorobjov, Nicolai N., 
{\it ``On the Number of Homotopy Types of Fibres of a Definable 
Map,''} Journal of the London Mathematical Society, to appear. 
Also available as Math ArXiV preprint {\tt math.AG/0605517} . 

\bibitem[BRS06]{thresh} Bihan, Frederic; Rojas, J.\ Maurice; and Stella, 
Casey, {\it ``First Steps in Algorithmic Fewnomial Theory,''} submitted for 
publication. Also available as Math ArXiV preprint {\tt math.AG/0411107} .

\bibitem[BS06]{bs} Bihan, Frederic and Sottile, Frank, {\it ``New 
Fewnomial Upper Bounds from Gale Dual \linebreak Polynomial Systems,''} 
Moscow Mathematical Journal, to appear. Also available as 
Math ArXiV preprint {\tt math.AG/0609544} .  

\bibitem[BCSS98]{bcss} Blum, Lenore; Cucker, Felipe; Shub, Mike; and
Smale, Steve, {\it Complexity and Real Computation,} Springer-Verlag, 1998. 

\bibitem[CSMP03]{pardo} Castro, David; San Mart\'\i{}n, Jorge; and Pardo, 
Luis-Miguel,  
{\it ``Systems of rational \linebreak polynomial equations have polynomial size
approximate zeros on the average,''}  J.\ Complexity  19  (2003), pp.\
161--209.

\bibitem[CG84]{chigo} Chistov, Alexander L., and Grigoriev, Dima Yu, {\it
``Complexity of Quantifier Elimination in the Theory of Algebraically
Closed Fields,''} Lect.\ Notes Comp.\ Sci.\ 176, Springer-Verlag (1984).

\bibitem[CZ02]{cohenzannier} Cohen, Paula B.\ and Zannier, Umberto, {\it 
``Fewnomials and intersections of lines with real analytic subgroups in 
$\mathbf{G}^n_m$,''}  
Bull.\ London Math.\ Soc.\ 34 (2002), no.\ 1, pp.\ 21--32.

\bibitem[Cos98]{coste} Coste, Michel, {\it ``Topological types of 
fewnomials,''} 
Singularities Symposium ---\L{}ojasiewicz 70 (Krak\'ow, 1996; Warsaw, 1996), 
pp.\ 81--92, Banach Center Publ., 44, Polish Acad.\ Sci., Warsaw, 1998. 

\bibitem[CD06]{cd} Cueto, Mar\'{\i}a Ang\'elica and Dickenstein, Alicia,
{\it ``Some results on inhomogeneous \linebreak discriminants,''} to appear: 
Proc.\ XVI 
CLA, Biblioteca de la Revista Matem\'atica Iberoamericana.

\bibitem[DFS05]{dfs} Dickenstein, Alicia; Feichtner, Eva Maria; 
and Sturmfels, Bernd, {\it ``Tropical \linebreak Discriminants,''} 
Math ArXiV preprint {\tt math.AG/0510126} . 

\bibitem[DS02]{codim2} Dickenstein, Alicia and Sturmfels, Bernd, 
{\sl ``Elimination Theory in Codimension Two,''} \linebreak Journal of 
Symbolic Computation (2002) {\bf 34}, pp.\ 119--135.

\bibitem[vdD98]{ominimal} van den Dries, Lou, {\it Tame topology and o-minimal 
structures,} London Mathematical \linebreak Society Lecture Note Series, 248, 
Cambridge University Press, Cambridge, 1998. 

\bibitem[Ful93]{tfulton} Fulton, William, {\it
Introduction to Toric Varieties}, Annals of Mathematics Studies, no.\ 131,
Princeton University Press, Princeton, New Jersey, 1993.

\bibitem[GV01]{gv} Gabrielov, Andrei and Vorobjov, Nicolai, 
{\it ``Complexity of cylindrical decompositions of sub-Pfaffian sets,''} 
Effective methods in algebraic geometry (Bath, 2000), 
J.\ Pure Appl.\ Algebra 164 (2001), no.\ 1--2, pp.\ 179--197.

\bibitem[GVZ04]{gabrielov} Gabrielov, Andrei; Vorobjov, Nikolai; and Zell, 
Thierry, {\it ``Betti Numbers of Semialgebraic and Sub-Pfaffian Sets,''} 
J.\ London Math.\ Soc.\ {\bf 69} (2004), pp.\ 27--43.  

\bibitem[GKZ94]{gkz94} Gel'fand, Israel Moseyevitch; Kapranov, Misha M.; and
Zelevinsky, Andrei V.;\linebreak 
{\it Discriminants, Resultants and Multidimensional Determinants,}
Birkh\"auser, Boston, 1994.

\bibitem[GK03]{rimas} {\it Algebraic Geometry and 
Geometric Modelling,} 
Proceedings of a conference in Vilnius, Lithuania (July 29-August 2, 2002), 
Contemporary Mathematics, (edited by Ron Goldman and Rimvydas Krasauskas), 
vol.\ 334, American Mathematical Society, 2003.

\bibitem[Haa02]{haas} Haas, Bertrand, {\it ``A Simple Counter-Example
to Kushnirenko's Conjecture,''} Beitr\"age zur Algebra und Geometrie,
Vol.\ 43, No.\ 1, pp.\ 1--8 (2002).

\bibitem[Hab48]{habicht} Habicht, Walter, {\it ``Eine Verallgemeinerung des
Sturmschen Wurzelzählverfahrens,''} \linebreak 
Comment.\ Math.\ Helv.\ {\bf 21} (1948), pp.\ 99--116.

\bibitem[Har76]{harnack} Harnack, Carl Gustav Axel, {\it ``\"Uber die 
Vielfaltigkeit der ebenen algebraischen Kurven,''} Math.\ Ann.\ 
{\bf 10} (1876), pp.\ 189--199. 

\bibitem[Kal03]{kaloshin} Kaloshin, V., {\it ``The existential Hilbert 16-th 
problem and an estimate for cyclicity of \linebreak elementary polycycles,''} 
Invent.\ Math.\ 151 (2003), no.\ 3, pp.\ 451--512.

\bibitem[Kap91]{kapranov} Kapranov, Misha, 
 {\sl ``A characterization of A-discriminantal hypersurfaces in terms of the
 logarithmic Gauss map,''} Mathematische Annalen, 290, 1991, pp.\ 277--285.

\bibitem[Kho80]{kho} Khovanskii, Askold G., {\it ``On a Class of
Systems of Transcendental Equations,''} Dokl.\
Akad.\ Nauk SSSR {\bf 255} (1980), no.\ 4, pp.\ 804--807;
English transl.\ in Soviet Math.\ Dokl.\ {\bf 22} (1980),
no.\ 3.

\bibitem[Kho91]{few} \underline{\hspace{\khov}}, {\it Fewnomials,}
AMS Press, Providence, Rhode Island, 1991.

\bibitem[Koi97]{koirandim} Koiran, Pascal, {\it ``Randomized and
Deterministic Algorithms for the Dimension of Algebraic Varieties,''}
Proceedings of the 38$^\thth$ Annual IEEE Computer Society
Conference on Foundations of Computer Science (FOCS),
Oct.\ 20--22, 1997, ACM Press.

\bibitem[Kri06]{krick} Krick, Teresa, {\it personal communication,} 
at the Institute for Mathematics and its \linebreak 
Applications (University of 
Minnesotta, Minneapolis), Sept.\ 15, 2006. 

\bibitem[LRW03]{tri} Li, Tien-Yien; Rojas, J.\ Maurice; and
Wang, Xiaoshen, {\it ``Counting Real Connected \linebreak 
Components of Trinomial
Curves Intersections and m-nomial Hypersurfaces,''} Discrete and
Computational Geometry, 30:379--414 (2003).

\bibitem[LM01]{lickroy} Lickteig, Thomas and Roy, Marie-Francoise,
{\it ``Sylvester-Habicht Sequences and Fast Cauchy Index Computation,''}
J.\ Symbolic Computation (2001) {\bf 31}, pp.\ 315--341.

\bibitem[Loe91]{loeser} Loeser, Fran\c{c}ois, {\it ``Polytopes secondaires 
et discriminants,"} S\'eminaire Bourbaki, Vol.\ 1990/91.  Ast\'erisque  
No.\ 201--203 (1991), Exp.\ No.\ 742, pp.\ 387--420 (1992). 

\bibitem[OK00]{orevkov} Orevkov, S.\ Yu.\ and Kharlamov, V.\ M., {\it 
``Asymptotic growth of the number of classes of real plane algebraic 
curves when the degree increases,''} J.\ of Math. Sciences,  113 (2003), no.\ 5,
pp.\ 666--674.

\bibitem[PT04]{pt} Passare, Mikael and Tsikh, August, {\it ``Algebraic 
equations and hypergeometric series,''} in: The Legacy of Niels Henrik Abel, 
Springer-Verlag, Berlin; Heidelberg, 2004.

\bibitem[PT05]{passare} Passare, Mikael and
Tsikh, August, {\it ``Amoebas: their spines and their contours,''} 
\linebreak 
Idempotent mathematics and mathematical physics, Contemp.\ Math., 
v.\ 377, Amer.\ Math.\ Soc., Providence, RI, 2005, pp.\ 275--288. 

\bibitem[Per05]{perrucci} Perrucci, Daniel, {\it ``Some Bounds for the 
Number of Components of Real Zero Sets of Sparse Polynomials,''} 
Discrete and Computational Geometry, to appear. 

\bibitem[Pla84]{plaisted} Plaisted, David A., {\it ``New NP-Hard and
NP-Complete Polynomial and Integer Divisibility Problems,''}
Theoret.\ Comput.\ Sci.\ 31 (1984), no.\ 1--2, 125--138.

\bibitem[Roj03]{why} Rojas, J.\ Maurice, {\it ``Why Polyhedra Matter in
Non-Linear Equation Solving,''} paper\linebreak 
corresponding to an invited talk delivered at a conference on Algebraic
Geometry and \linebreak Geometric Modelling (Vilnius, Lithuania, 
July 29 -- August 2, 2002),
Contemporary Mathematics, vol.\ 334, pp.\ 293--320, AMS Press, 2003.

\bibitem[Roj04]{amd} Rojas, J.\ Maurice, {\it ``Arithmetic Multivariate
Descartes' Rule,''} American Journal of \linebreak 
Mathematics, vol.\ 126, no.\ 1, February 2004, pp.\ 1--30.

\bibitem[RY05]{rojasye} Rojas, J.\ Maurice and Ye, Yinyu, {\it
``On Solving Sparse Polynomials in Logarithmic Time,''}
Journal of Complexity, special issue for the
2002 Foundations of Computation Mathematics (FOCM) meeting,
February 2005, pp.\ 87--110.

\bibitem[RS98]{pole} Rosenthal, Joachim and Sottile, Frank, {\it 
``Some remarks on real and complex output feedback,''} 
Systems Control Lett.\ 33 (1998), no.\ 2, pp.\ 73--80. 

\bibitem[Roy96]{marie} Roy, Marie-Fran\c{c}oise, {\it ``Basic Algorithms in
Real Algebraic Geometry and their Complexity: from Sturm's Theorem to the
Existential Theory of Reals,''} Lectures in Real Geometry (Madrid, 1994), pp.\
1--67, de Gruyter Exp.\ Math., 23, de Gruyter, Berlin, 1996.

\bibitem[Sma86]{smale} Smale, Steve, {\it
``Newton's Method Estimates from Data at One Point,''}
The Merging of \linebreak Disciplines: New Directions in Pure, Applied, and
Computational Mathematics (Laramie, Wyo., 1985), pp.\ 185--196, Springer, New
York, 1986.

\bibitem[Sma00]{21} \underline{\hspace{\smale}}, {\it
``Mathematical Problems for the Next Century,''}
Mathematics: Frontiers and Perspectives, pp.\ 271--294, Amer.\ Math.\ Soc.,
Providence, RI, 2000.

\bibitem[SL54]{descartes} Smith, David Eugene and Latham, Marcia L., {\it
The Geometry of Ren\'e Descartes,} \linebreak 
translated from the French and Latin
(with a facsimile of Descartes' 1637 French edition),
Dover Publications Inc., New York (1954).

\bibitem[Sto98]{storjo} Storjohann, Arne, {\it ``Computing Hermite and 
Smith normal forms of triangular integer \linebreak matrices,''} Linear Algebra 
Appl.\ 282 (1998), no.\ 1--3, pp.\ 25--45. 

\bibitem[Stu35]{sturm} Sturm, C., {\it
``M\'emoire sur la r\'esolution des \'equations num\'eriques,''}
Inst.\ France Sc.\ Math.\ Phys., {\bf 6} (1835).

\bibitem[VG03]{butterfly} Vakulenko, Sergey and Grigoriev, Dmitry, 
{\it ``Complexity of gene circuits, Pfaffian functions and the 
morphogenesis problem,''} 
C.\ R.\ Math.\ Acad.\ Sci.\ Paris 337 (2003), no.\ 11, pp.\ 721--724. 

\bibitem[Ver06]{jan} Verschelde, Jan, {\it PHC-Pack,} 
{\tt http://www.math.uic.edu/\~{}jan/download.html } . 

\bibitem[Wil99]{wilkie} Wilkie, A.\ J., {\it ``A theorem of the complement 
and some new o-minimal structures,''} Selecta Math.\ (N.S.)  5  (1999),  
no.\ 4, pp.\ 397--421.

\end{thebibliography}

\end{document}